\documentclass[11pt]{article}
\usepackage{amssymb}

\newtheorem{theorem}{{\bf Theorem}}[section]
\newtheorem{lemma}{{\bf Lemma}}[section]

\newtheorem{definition}{{\bf Definition}}[section]
\newtheorem{corollary}{{\bf Corollary}}[section]
\newtheorem{remark}{{\bf Remark}}[section]

\begin{document}

\title{Reflected backward SDEs with two barriers under monotonicity and general
increasing conditions}
\author{Mingyu Xu \thanks{%
Email: xvmingyu@gmail.com} \\
%EndAName
{\small Departement des Math\'ematiques, Universit\'e du Maine, Le Mans
France;}\\
{\small Department of Financial Mathematics and Control science, School of
Mathematical Science} \\
{\small Fudan University, Shanghai, 200433, China.}}
\date{}
\maketitle

\textbf{Abstract} In this paper, we prove the existence and uniqueness
result of the reflected BSDE with two continuous barriers under monotonicity
and general increasing condition on $y$, with Lipschitz condition on $z$.


\section{Introduction}

Nonlinear backward stochastic differential equations (BSDEs in short) were
firstly introduced by Pardoux and Peng (1990), who proved the existence and
uniqueness of adapted solutions, when the coefficient $f$ is Lipschitz in $%
(y,z)$ uniformly in $(t,\omega )$, with square-integrability assumptions on
the coefficient $f(t,\omega ,y,z)$ and terminal condition $\xi $. Later,
Pardoux (1999) and Briand, Delyon, Hu, Pardoux and Stoica (2003) studied the
solution of a BSDE with a coefficient $f(t,\omega ,y,z)$ that satisfies only
monotonicity, continuity and general increasing growth conditions with
respect to $y$, and Lipschitz on $z$. That is, for some real number $\mu \in
\mathbb{R}$, $k\geq 0$ and some continuous increasing function $\varphi :%
\mathbb{R}_{+}\rightarrow \mathbb{R}_{+}$: $\forall t\in [0,T]$, $%
y,y^{\prime }\in \mathbb{R}$, $z,z^{\prime }\in \mathbb{R}^{d}$:
\begin{eqnarray}
\left| f(t,y,z)\right| &\leq &\left| f(t,0,z)\right| +\varphi (\left|
y\right| ),  \label{intro} \\
(y-y^{\prime })(f(t,y,z)-f(t,y^{\prime },z)) &\leq &\mu \left| y-y^{\prime
}\right| ^{2},  \nonumber \\
\left| f(t,y,z)-f(t,y,z^{\prime })\right| &\leq &k\left| z-z^{\prime
}\right| .  \nonumber
\end{eqnarray}

Reflected backward stochastic differential equations (RBSDEs in short) with
one lower barrier were studied by El Karoui, Kapoudjian, Pardoux, Peng and
Quenez (1997), in one dimension. The solution is constrained to remain above
a continuous lower-boundary process with the help of an continuous
increasing process. Later, Cvitanic and Karatzas (1996) studied the backward
stochastic differential equation with two barriers. A solution to such
equation associated to a terminal condition $\xi $, a coefficient $%
f(t,\omega ,y,z)$ and two barriers $L$ and $U$, is a triple
$(Y,Z,K)$ of adapted processes, valued in $\mathbb{R}^{1+d+1}$,
which satisfies $\;$

\[
Y_{t}=\xi
+\int_{t}^{T}f(s,Y_{s},Z_{s})ds+K_{T}^{+}-K_{t}^{+}-(K_{T}^{-}-K_{t}^{-})-%
\int_{t}^{T}Z_{s}dB_{s},\;0\leq t\leq T\mbox{ a.s}.
\]
$L_{t}\leq Y_{t}\leq U_{t}$, $0\leq t\leq T$ and $%
\int_{0}^{T}(Y_{s}-L_{s})dK_{s}^{+}=\int_{0}^{T}(Y_{s}-U_{s})dK_{s}^{-}=0,$
a.s. In this case, a solution $Y$ has to remain between the lower boundary $%
L $ and upper boundary $U$, almost surely. This is achieved by the
cumulative action of two continuous, increasing reflecting processes $K^{\pm
}$, which act in a minimal way when $Y$ attempts to cross barriers. And the
authors proved the existence and uniqueness of the solution when $f(t,\omega
,y,z)$ is Lipschitz on $(y,z)$ uniformly in $(t,\omega )$ and when $L<U$ on $%
[0,T]$ and there exists a different supermartingale between $L$ and $U$
(Mokobodski's assumption in Dynkin game). Furthermore they established the
connection between solution $Y$ and the value of Dynkin games (certain
stochastic games of stopping). Then in \cite{HLM}, the existence of a
solution was proved when $f$ is only continuous with linear growth in $(y,z)$%
, but in the case when one obstacle is smooth. Later, Lepeltier and San
Martin used the penalization method to prove the existence of a solution to
such equation, with same assumption on $f$ as in \cite{HLM}, without extra
smoothness of the barriers, i.e. when $L$ and $U$ are continuous, $L<U$ on $%
[0,T]$, and Mokobodski's assumption.

More recently, Lepeltier, Matoussi and Xu proved the existence and
uniqueness of the solution to the reflected BSDE with one lower continuous
barrier under the assumption (\ref{intro}) for $f$. The existence is proved
by approximation. In this paper, we consider the reflected BSDE with two
continuous barrier under the assumption (\ref{intro}), and give the
uniqueness and existence of the solution, which is obtained by approximation.

The paper is organized as following: In subsection 2.1, we present
notations and assumptions; then we prove the main results of this
paper, the existence and uniqueness of the solution in subsection
2.2; in subsection 2.3 we prove an important theorem for the
existence in five steps. Finally, in section 3, we prove several
comparison theorems with respect to RBSDE with one or two barriers,
which are used in the proof of existence.

\section{RBSDE's with two continuous barriers}

\subsection{Assumptions and notations}

Let $(\Omega ,\mathcal{F},P)$ be a complete probability space, and $%
B=(B_{1},B_{2},\cdots ,B_{d})^{\prime }$ be a $d$-dimensional Brownian
motion defined on a finite interval $[0,T]$, $0<T<+\infty $. Denote by $\{%
\mathcal{F}_{t};0\leq t\leq T\}$ the natural filtration generated by the
Brownian motion $B:$%
\[
\mathcal{F}_{t}=\sigma \{B_{s};0\leq s\leq t\},
\]
where $\mathcal{F}_{0}$ contains all $P-$null sets of $\mathcal{F}$.

We denote the following notations. For any given $n\in \mathbf{N}$, let us
introduce the following spaces:
\[
\begin{array}{ll}
\mathbf{L}_{n}^{2}(\mathcal{F}_{t})= & \{\xi :n\mbox{-dimensional }\mathcal{F%
}_{t}\mbox{-measurable random variable, s.t. }E(|\xi |^{2})<+\infty \}, \\
\mathbf{H}_{n}^{2}(0,T)= & \{\psi :n\mbox{-dimensional }\mathcal{F}_{t}\mbox{%
-predictable process on the interval }[0,T]\mbox{, } \\
& \mbox{s.t. }E\int_{0}^{T}\left\| \psi (t)\right\| ^{2}dt<+\infty \}, \\
\mathbf{S}_{n}^{2}(0,T)= & \{\psi :n\mbox{-dimensional }\mathcal{F}_{t}\mbox{%
-progressively measurable continuous process } \\
& \mbox{on the interval }[0,T]\mbox{, s.t. }E(\sup_{0\leq t\leq
T}\left\|
\psi (t)\right\| ^{2})<+\infty \}, \\
\mathbf{A}^{2}(0,T)= & \{K:\mbox{real valued
}\mathcal{F}_{t}\mbox{-adapted
increasing continuous process, s.t. }K(0)=0\mbox{,} \\
& \mbox{and }E(K(T)^{2})<+\infty \}. \\
\mathbf{VF}^{2}(0,T)= & \{V:\mbox{real valued
}\mathcal{F}_{t}\mbox{-adapted
continuous process with finite variation, s.t. } \\
& V=K^{+}-K^{-}\mbox{, with }K^{\pm }\in \mathbf{A}^{2}(0,T)\}.
\end{array}
\]

Finally, we shall denote by $\mathcal{P}$ the $\sigma $-algebra of
predictable sets on $[0,T]\times \Omega $. In the real--valued case, i.e., $%
n=1$, these spaces will be simply denoted by $\mathbf{L}^{2}(\mathcal{F}%
_{t}) $, $\mathbf{H}^{2}(0,T)$ and $\mathbf{S}^{2}(0,T)$, respectively.

Let us consider the reflected backward stochastic differential equation with
monotonic condition in $y$ on a fixed time interval; we need the following
assumptions:

\noindent\textbf{Assumption 2.1.} A final condition $\xi \in \mathbf{L}^{2}(\mathcal{F%
}_{T})$.

\noindent\textbf{Assumption 2.2.} A coefficient $f:\Omega \times [0,T]\times \mathbf{%
R\times R}^{d}\rightarrow \mathbb{R}$, satisfying for some
continuous increasing function $\varphi :\mathbb{R}_{+}\rightarrow
\mathbb{R}_{+}$, real numbers $\mu \in \mathbb{R}$ and $k>0$:

$
\begin{array}{ll}
\mbox{(i)} & f(\cdot ,y,z)\mbox{ is progressively measurable,
}\forall
(y,z)\in \mathbb{R\times R}^{d}; \\
\mbox{(ii)} & \left| f(t,y,z)\right| \leq \left| f(t,0,z)\right|
+\varphi
(\left| y\right| )\mbox{, }\forall (t,y,z)\in [0,T]\times \mathbb{R\times R}%
^{d}\mbox{, a.s.;} \\
\mbox{(iii)} & E\int_{0}^{T}\left| f(t,0,0)\right| ^{2}dt<\infty ; \\
\mbox{(iv)} & \left| f(t,y,z)-f(t,y,z^{\prime })\right| \leq k\left|
z-z^{\prime }\right| \mbox{, }\forall (t,y)\in [0,T]\times
\mathbb{R}\mbox{,
}z,z^{\prime }\in \mathbb{R}^{d}\mbox{, a.s.} \\
\mbox{(v)} & (y-y^{\prime })(f(t,y,z)-f(t,y^{\prime },z))\leq \mu
\left|
y-y^{\prime }\right| ^{2}\mbox{, }\forall (t,z)\in [0,T]\times \mathbb{R}^{d}%
\mbox{, }y,y^{\prime }\in \mathbb{R}\mbox{, a.s.} \\
\mbox{(vi)} & y\rightarrow f(t,y,z)\mbox{ is continuous, }\forall
(t,z)\in [0,T]\times \mathbb{R}^{d}\mbox{, a.s.}
\end{array}
$

\noindent\textbf{Assumption 2.3.} Two barriers $L_{t}$, $U_{t}$, which are $\mathcal{F%
}_{t}$-progressively measurable continuous processes, defined on the
interval $[0,T]$, satisfying

(i)

\[
E[\varphi ^{2}(\sup_{0\leq t\leq T}(e^{\mu t}(L_{t})^{+}))]<\infty
,E[\varphi ^{2}(\sup_{0\leq t\leq T}(e^{\mu t}(U_{t})^{-}))]<\infty ,
\]
$(L)^{+},(U)^{-}\in \mathbf{S}^{2}(0,T)$, and $L_{T}\leq \xi \leq U_{T}$,
a.s., where $(L)^{+}$(resp. $(U)^{-}$) is the positive part (resp. negative)
part of $L$ (resp. $U$).

(ii) there exists a process $J_t=J_0+\int_0^t\phi _sdB_s-V_t^{+}+V_t^{-}$, $%
J_T=\xi $ with $\phi \in \mathbf{H}_d^2(0,T)$, $V^{+},V^{-}\in \mathbf{A}%
^2(0,T)$, s.t.
\[
L_t\leq J_t\leq U_t\mbox{, for }0\leq t\leq T.
\]

(iii) $L_{t}<U_{t}$, a.s., for $0\leq t<T.$

Now we introduce the definition of the solution of RBSDE with two barriers $%
L $ and $U$\textbf{.}

\begin{definition}
\label{Def-tb}We say that $(Y_{t},Z_{t},K_{t})_{0\leq t\leq T}$ is a
solution of the backward stochastic differential equation with two
continuous reflecting barriers $L(\cdot )$ and $U(\cdot )$, terminal
condition $\xi $ and coefficient $f$, which is denoted as RBSDE$(\xi ,f,L,U)$%
, if the followings hold:

(1) $Y\in \mathbf{S}^{2}(0,T)$, $Z\in \mathbf{H}_{d}^{2}(0,T)$, and $K\in
\mathbf{VF}^{2}(0,T)$, $K=K^{+}-K^{-}$, where $K^{\pm }\in \mathbf{A}%
^{2}(0,T)$.

(2) $Y_{t}=\xi
+\int_{t}^{T}f(s,Y_{s},Z_{s})ds+K_{T}^{+}-K_{t}^{+}-(K_{T}^{-}-K_{t}^{-})-%
\int_{t}^{T}Z_{s}dB_{s},\;\;0\leq t\leq T$ a.s.

(3) $L_{t}\leq Y_{t}\leq U_{t},\;\;0\leq t\leq T,$ a.s.

(4) $\int_{0}^{T}(Y_{s}-L_{s})dK_{s}^{+}=%
\int_{0}^{T}(Y_{s}-U_{s})dK_{s}^{-}=0,$ a.s.
\end{definition}

Actually, a general solution of our RBSDE$(\xi ,f,L,U)$ would satisfy the
assumptions (1) to (4). The state-process $Y(\cdot )$ is forced to stay
between the barrier $L(\cdot )$ and $U(\cdot )$, thanks to the cumulation
action of the reflection processes $K^{+}(\cdot )$ and $K^{-}(\cdot )$
respectively, which act only when necessary to prevent $Y(\cdot )$ from
crossing the respective barrier, and in this sense, its action can be
considered minimal, i.e. the integrability assumption (4). From the fact
that $K^{\pm }\in \mathbf{A}^{2}(0,T)$ is continuous and (2), it follows
that $Y$ is continuous.

\begin{remark}
We have an analogue result of Proposition 4.1 in \cite{CK}. Precisely, the
square-integrable solution $Y$ of the RBSDE$(\xi ,f,L,U)$ is the value of
the Dynkin game problem, whose payoff is
\[
R_{t}(\sigma ,\tau )=\int_{t}^{\sigma \wedge \tau }f(s,Y_{s},Z_{s})ds+\xi
1_{\{\sigma \wedge \tau =T\}}+L_{\tau }1_{\{\tau <T,\tau \leq \sigma
\}}+U_{\sigma }1_{\{\sigma <\tau \}},
\]
and a saddle-point $(\widehat{\sigma }_{t},\widehat{\tau }_{t})\in \mathcal{T%
}_{t}\times \mathcal{T}_{t}$ is given by
\begin{eqnarray*}
\widehat{\sigma }_{t} &=&\inf \{s\in [t,T);Y_{s}=U_{s}\}\wedge T, \\
\widehat{\tau }_{t} &=&\inf \{s\in [t,T);Y_{s}=L_{s}\}\wedge T.
\end{eqnarray*}
\end{remark}

\subsection{Main results}

Our main results in this paper is following:

\begin{theorem}
\label{unique}Under the assumptions 2.1, 2.2 and 2.3, the RBSDE$(\xi ,f,L,U)$
has the unique solution $(Y,Z,K)$, which satisfies definition \ref{Def-tb}
(1)-(4).
\end{theorem}

\textbf{Proof. }\textsl{Uniqueness. }Suppose that the triples $(Y,Z,K)$ and $%
(Y^{\prime },Z^{\prime },K^{\prime })$ are two solutions of the RBSDE$(\xi
,f,L)$, i.e. satisfy (1)-(4) of definition \ref{Def-tb}. Set $\Delta
Y=Y-Y^{\prime }$, $\Delta Z=Z-Z^{\prime }$, $\Delta K=\Delta K-\Delta
K^{\prime }$, with $\Delta K^{+}=K^{+}-K^{+\prime }$, $\Delta
K^{-}=K^{-}-K^{-\prime }$. Applying It\^{o}'s formula to $\Delta Y^{2}$ on
the interval $[t,T]$, and taking expectation on both sides, it follows
\[
E\left| \Delta Y_{t}\right| ^{2}+E\int_{t}^{T}\left| \Delta Z_{s}\right|
^{2}ds\leq 2(k^{2}+\mu )E\int_{t}^{T}\Delta Y_{s}^{2}ds+\frac{1}{2}%
E\int_{t}^{T}\left| \Delta Z_{s}\right| ^{2}ds,
\]
in view of monotonic assumption on $y$, Lipschitz assumption on $z$, and $%
\int_{t}^{T}\Delta Y_{s}d\Delta K_{s}\leq 0$. We get
\[
E\left| \Delta Y_{t}\right| ^{2}\leq 2(k^{2}+\mu )E\int_{t}^{T}\Delta
Y_{s}^{2}ds.
\]
From the Gronwall's inequality, it follows $E\left| \Delta Y_{t}\right|
^{2}=E\left| Y_{t}-Y_{t}^{\prime }\right| ^{2}=0$, $0\leq t\leq T$, i.e. $%
Y_{t}=Y_{t}^{\prime }$ a.s.; then we have also $E\int_{0}^{T}\left| \Delta
Z_{s}\right| ^{2}ds=E\int_{0}^{T}\left| Z_{s}-Z_{s}^{\prime }\right|
^{2}ds=0 $, from which follows $K_{t}=K_{t}^{\prime }$.

\textsl{Existence. }We firstly present the following existence theorem when $%
f$ does not depend on $z$, which will be proved a little later.

\begin{theorem}
\label{exist22}Suppose that $\xi $, $f$ and $L$, $U$ satisfy assumption 2.1,
2.2 and 2.3, then for any process $Q\in \mathbf{H}_{d}^{2}(0,T)$, there
exists a unique triple of progressively measurable processes $%
\{(Y_{t},Z_{t},K_{t})_{0\leq t\leq T}\}\in \mathbf{S}^{2}(0,T)\times \mathbf{%
H}_{d}^{2}(0,T)\times \mathbf{VF}^{2}(0,T)$, with $K=K^{+}-K^{-}$, $%
(K_{t}^{\pm })_{0\leq t\leq T}\in \mathbf{A}^{2}(0,T)$, which satisfies \ref
{Def-tb} (1), (3), (4) and
\[
Y_{t}=\xi
+\int_{t}^{T}f(s,Y_{s},Q_{s})ds+K_{T}^{+}-K_{t}^{+}-(K_{T}^{-}-K_{t}^{-})-%
\int_{t}^{T}Z_{s}dB_{s},0\leq t\leq T\mbox{.}
\]
\end{theorem}

Thanks to Theorem \ref{exist22}, we can construct a mapping $\Phi $ from $%
\mathcal{S}$ into itself, where $\mathcal{S}$ is defined as the
space of the progressively measurable processes
$\{(Y_{t},Z_{t});0\leq t\leq T\},$ valued in $\mathbb{R\times
R}^{d}$ which satisfy (1) as follows.

Given $(P,Q)\in \mathcal{S}$, $(Y,Z)=\Phi (P,Q)$ is the unique solution of
following RBSDE
\[
Y_{t}=\xi
+\int_{t}^{T}f(s,Y_{s},Q_{s})ds+K_{T}-K_{t}-\int_{t}^{T}Z_{s}dB_{s},
\]
i.e., if we define the process
\[
K_{t}=Y_{t}-Y_{0}-\int_{0}^{t}f(s,Y_{s},Q_{s})ds+\int_{0}^{t}Z_{s}dB_{s},0%
\leq t\leq T,
\]
then the triple $(Y,Z,K)$ satisfies definition \ref{Def-tb} (1)-(4), with $%
f(s,y,z)=f(s,y,Q_{s})$.

Consider another element of $\mathcal{S}$, and define $(Y^{\prime
},Z^{\prime })=\Phi (P^{\prime },Q^{\prime })$; set
\begin{eqnarray*}
\Delta P &=&P-P^{\prime },\Delta Q=Q-Q^{\prime },\Delta Y=Y-Y^{\prime
},\Delta Z=Z-Z^{\prime }, \\
\Delta K &=&K^{+}-K^{-},\Delta K^{+}=K^{+}-K^{+\prime },\Delta
K^{-}=K^{-}-K^{-\prime }.
\end{eqnarray*}
We apply the It\^{o}'s formula to $e^{\gamma t}\left| \Delta Y_{t}\right|
^{2}$ on the interval $[t,T]$, for $\gamma >0$,
\begin{eqnarray*}
&&e^{\gamma t}E\left| \Delta Y_{t}\right| ^{2}+E\int_{t}^{T}e^{\gamma
s}(\gamma \left| \Delta Y_{s}\right| ^{2}+\left| \Delta Z_{s}\right| ^{2})ds
\\
&\leq &2(k^{2}+\mu )E\int_{t}^{T}e^{\gamma s}\left| \Delta Y_{s}\right|
^{2}ds+\frac{1}{2}E\int_{t}^{T}e^{\gamma s}\left| \Delta Q_{s}\right| ^{2}ds,
\end{eqnarray*}
since $\int_{t}^{T}e^{\gamma s}\Delta Y_{s}d\Delta
K_{s}=\int_{t}^{T}e^{\gamma s}\Delta Y_{s}d\Delta
K_{s}^{+}-\int_{t}^{T}e^{\gamma s}\Delta Y_{s}d\Delta K_{s}^{-}\leq 0$.
Hence, if we choose $\gamma =1+2(k^{2}+\mu )$, it follows
\begin{eqnarray*}
E\int_{t}^{T}e^{\gamma s}(\left| \Delta Y_{s}\right| ^{2}+\left| \Delta
Z_{s}\right| ^{2})ds &\leq &\frac{1}{2}E\int_{t}^{T}e^{\gamma s}\left|
\Delta Q_{s}\right| ^{2}ds \\
&\leq &\frac{1}{2}E\int_{t}^{T}e^{\gamma s}(\left| \Delta P_{s}\right|
^{2}+\left| \Delta Q_{s}\right| ^{2})ds.
\end{eqnarray*}
Consequently, $\Phi $ is a strict contraction on $\mathcal{S}$ equipped with
the norm
\[
\left\| (Y,Z)\right\| _{\gamma }=\left[ E\int_{0}^{T}e^{\gamma s}(\left|
Y_{s}\right| ^{2}+\left| Z_{s}\right| ^{2})ds\right] ^{\frac{1}{2}},
\]
and has a fixed point, which is the unique solution of the RBSDE$(\xi
,f,L,U) $. $\square $

\subsection{Proof of theorem \ref{exist22}}

Now we prove the theorem \ref{exist22} in several steps for the existence of
solution. We write $f(s,y)$ for $f(s,y,Q_{s})$. First we note that the
triple $(Y,Z,K)$ solves the RBSDE$(\xi ,f,L,U)$, $K=K^{+}-K^{-}$, if and
only if
\begin{equation}
(\overline{Y}_{t},\overline{Z}_{t},\overline{K}_{t}^{+},\overline{K}%
_{t}^{-}):=(e^{\lambda t}Y_{t},e^{\lambda t}Z_{t},\int_{0}^{t}e^{\lambda
s}dK_{s}^{+},\int_{0}^{t}e^{\lambda s}dK_{s}^{-})  \label{trans}
\end{equation}
solves the RBSDE$(\overline{\xi },\overline{f},\overline{L},\overline{U})$,
where
\[
(\overline{\xi },\overline{f}(t,y),\overline{L}_{t},\overline{U}_{t})=(\xi
e^{\lambda T},e^{\lambda t}f(t,e^{-\lambda t}y)-\lambda y,e^{\lambda
t}L_{t},e^{\lambda t}U_{t}).
\]

If we choose $\lambda =\mu $, then the coefficient $\overline{f}$ satisfies
the same assumptions in assumption 2.2 as $f$, but with assumption 2.2-(v)
replaced by

(v') $(y-y^{\prime })(f(t,y,z)-f(t,y^{\prime },z))\leq 0$.

Since we are in $1$-dimensional case, (v') means that $f$ is decreasing on $%
y $. From another part the barriers $\overline{L}$, $\overline{U}$ satisfies:

\textbf{(i'):}
\begin{eqnarray*}
E[\sup_{0\leq t\leq T}(\overline{L}_{t})^{+}] &<&\infty ,E[\varphi
^{2}(\sup_{0\leq t\leq T}(\overline{L}_{t})^{+})]=E[\varphi ^{2}(\sup_{0\leq
t\leq T}(e^{\mu t}(L_{t})^{+}))]<\infty , \\
E[\sup_{0\leq t\leq T}(\overline{U}_{t})^{-}] &<&\infty ,E[\varphi
^{2}(\sup_{0\leq t\leq T}(\overline{U}_{t})^{-})]=E[\varphi ^{2}(\sup_{0\leq
t\leq T}(e^{\mu t}(U_{t})^{-}))]<\infty .
\end{eqnarray*}

In the following, we shall work with assumption 2.2'\textbf{\ }which is
assumption 2.2 with (v) replaced by (v') and assumption 2.3' which is
assumption 2.2 with \textbf{(i') }instead of \textbf{(i)}.

\textbf{Proof of Theorem \ref{exist22}: }First, let us recall the
assumptions on the coefficient $f$:

\noindent\textbf{Assumption 2.4. }For $y\in \mathbb{R}$, $s\in
[0,T]$,
\[
\begin{array}{ll}
\mbox{(i)} & \left| f(s,y)\right| \leq \left| f(s,0,0)\right|
+k\left|
Q_{s}\right| +\varphi (\left| y\right| ); \\
\mbox{(ii)} & E\int_{0}^{T}\left| f(t,0)\right| ^{2}dt<\infty ; \\
\mbox{(iii)} & (y-y^{\prime })(f(s,y)-f(s,y^{\prime }))\leq 0; \\
\mbox{(iv)} & y\rightarrow f(s,y)\mbox{ is continuous, a.s..}
\end{array}
\]

We point out that we always denote by $c>0$ a constant whose value can be
changed line by line. The proof will be done by five steps as following.

\begin{itemize}
\item  Using a penalization method we prove the existence under the
assumption
\begin{equation}
\left| \xi \right| +\sup_{0\leq t\leq T}\left| f(t,0)\right| +\sup_{0\leq
t\leq T}L_{t}^{+}+\sup_{0\leq t\leq T}U_{t}^{-}\leq c.  \label{assup-bon}
\end{equation}

\item  Approximating the lower barrier $L$, we prove the existence under the
assumption that $L$ satisfies assumption 2.3'\textbf{-(i)} and the bounded
assumption on $\xi $, $f(t,0)$ and $\sup_{0\leq t\leq T}U_{t}^{-}$.

\item  Like above step, we approximate the upper barrier $U$ to prove the
existence under assumption 2.3' and $\xi $ and $f(t,0)$ satisfy
\begin{equation}
\left| \xi \right| ^{2}+\sup_{0\leq t\leq T}\left| f(t,0)\right| ^{2}\leq c.
\label{assup-ini}
\end{equation}

\item  By approximation, we prove the existence of the solution under the
assumption $\xi \geq c$, $\inf_{0\leq t\leq T}f(t,0)\geq c$.

\item  Finally, we prove the existence of the solution under the assumption $%
\xi \in \mathbf{L}^{2}(\mathcal{F}_{T})$, $f(t,0)\in \mathbf{H}^{2}(0,T)$,
by approximation.
\end{itemize}

In each step, we use monotonic property of approximation solutions to get
the convergence.

\paragraph{\textbf{Step 1.}}

Consider the penalization equations with respect to the two barriers $L$, $U$%
, for $m$, $n\in \mathbf{N}$,
\begin{equation}
Y_{t}^{m,n}=\xi
+\int_{t}^{T}f(s,Y_{s}^{m,n})ds+m\int_{t}^{T}(Y_{s}^{m,n}-L_{s})^{-}ds-n%
\int_{t}^{T}(U_{s}-Y_{s}^{m,n})^{-}ds-\int_{t}^{T}Z_{s}^{m,n}dB_{s}.
\label{pen-two1}
\end{equation}
Set $f_{m,n}(s,y)=f(s,y)+m(y-L_{s})^{-}-n(U_{s}-y)^{-}$, obviously, $f_{m,n}$
satisfies the condition of Proposition 2.4 in \cite{Pardoux99}. So by the
Proposition 2.4 in \cite{Pardoux99}, there exists $%
(Y_{t}^{m,n},Z_{t}^{m,n})_{0\leq t\leq T}$, which is the solution of (\ref
{pen-two1}). Denote $K_{t}^{m,n,+}=m\int_{0}^{t}(Y_{s}^{m,n}-L_{s})^{-}ds$, $%
K_{t}^{m,n,-}=n\int_{0}^{t}(U_{s}-Y_{s}^{m,n})^{-}ds$.

Now let us do the uniformly a priori estimation of $%
(Y^{m,n},Z^{m,n},K^{m,n,+},K^{m,n,-})$.

\begin{lemma}
There exists a constant $C_{0}$ independent of $n$, such that
\[
E[\sup_{0\leq t\leq T}\left| Y_{t}^{m,n}\right| ^{2}+\int_{0}^{T}\left|
Z_{s}^{m,n}\right| ^{2}ds+(K_{T}^{m,n,+})^{2}+(K_{T}^{m,n,-})^{2}]\leq C_{0}.
\]
\end{lemma}

\textbf{Proof. }Consider the RBSDE$(\xi ,f,L)$ with one lower barrier $L$;
due to theorem 2.3 in \cite{LMX}, it admits a unique solution $(\overline{Y}%
_{t},\overline{Z}_{t},\overline{K}_{t})_{0\leq t\leq T}\in \mathbf{S}%
^{2}(0,T)\times \mathbf{H}_{d}^{2}(0,T)\times \mathbf{A}^{2}(0,T)$, which
satisfies
\begin{equation}
\overline{Y}_{t}=\xi +\int_{t}^{T}f(s,\overline{Y}_{s})ds+\overline{K}_{T}-%
\overline{K}_{t}-\int_{t}^{T}\overline{Z}_{s}dB_{s},  \label{ref-low1}
\end{equation}
$\overline{Y}_{t}\geq L_{t}$, $0\leq t\leq T$, $\int_{0}^{T}(\overline{Y}%
_{s}-L_{s})d\overline{K}_{s}=0$. In order to compare (\ref{ref-low1}) and (%
\ref{pen-two1}), we consider the penalization equation associated with the
RBSDE (\ref{ref-low1}), for $m\in \mathbf{N}$,
\begin{equation}
\overline{Y}_{t}^{m}=\xi +\int_{t}^{T}f(s,\overline{Y}_{s}^{m})ds+m%
\int_{t}^{T}(L_{s}-\overline{Y}_{s}^{m})^{+}ds-\int_{t}^{T}\overline{Z}%
_{s}^{m}dB_{s}.  \label{ref-low2}
\end{equation}
Comparing (\ref{pen-two1}) and (\ref{ref-low2}), we get $Y_{t}^{m,n}\leq
\overline{Y}_{t}^{m}$, $\forall t\in [0,T]$, $n\in \mathbf{N}$. Thank to the
convergence result of step1 and step 2 in the proof of theorem 2.3 in \cite
{LMX}, i.e. $\overline{Y}^{m}\rightarrow \overline{Y}$ in $\mathbf{S}%
^{2}(0,T)$. So we get$\ $for any $m,n\in \mathbf{N}$, $t\in [0,T]$, $%
Y_{t}^{m,n}\leq \overline{Y}_{t}$.

Similarly, we consider the RBSDE$(\xi ,f,U)$ with one upper barrier $U$.
There exists $(\underline{Y}_{t},\underline{Z}_{t},\underline{K}_{t})_{0\leq
t\leq T}\in \mathbf{S}^{2}(0,T)\times \mathbf{H}_{d}^{2}(0,T)\times \mathbf{A%
}^{2}(0,T)$, which satisfies
\begin{equation}
\underline{Y}_{t}=\xi +\int_{t}^{T}f(s,\underline{Y}_{s})ds-(\underline{K}%
_{T}-\underline{K}_{t})-\int_{t}^{T}\underline{Z}_{s}dB_{s},  \label{ref-up1}
\end{equation}
$\underline{Y}_{t}\leq U_{t}$, $0\leq t\leq T$, $\int_{0}^{T}(\underline{Y}%
_{s}-U_{s})d\underline{K}_{s}=0$. By the penalization equation associated
with (\ref{ref-up1}) and the comparison theorem, we deduce that $%
Y_{t}^{m,n}\geq \underline{Y}_{t}$, for any $m,n\in \mathbf{N}$, $t\in [0,T]$%
. Then we get, with the results of the step 1 in the proof of theorem 2.3
\cite{LMX},
\begin{equation}
\sup_{0\leq t\leq T}\left| Y_{t}^{m,n}\right| \leq \max \{\sup_{0\leq t\leq
T}\left| \overline{Y}_{t}\right| ,\sup_{0\leq t\leq T}\left| \underline{Y}%
_{t}\right| \}\leq C.  \label{est-mn}
\end{equation}
In the following, notice that assumption 2.4-(iii) implies that $f$ is
decreasing on $y$, for $s\in [0,T]$, so $f(s,\underline{Y}_{s})\geq
f(s,Y_{s}^{m,n})\geq f(s,\overline{Y}_{s})$, with the square-integrable
results of (\ref{ref-low1}) and (\ref{ref-up1}), it follows
\begin{equation}
\left| f(s,Y_{s}^{m,n})\right| \leq \max \{\left| f(s,\overline{Y}%
_{s})\right| ,\left| f(s,\underline{Y}_{s})\right| \}\leq C.  \label{est-fmn}
\end{equation}

To get the estimation of $(K^{m,n,+},K^{m,n,-},Z^{m,n})$, we apply It\^{o}'s
formula to $\left( Y^{m,n}\right) ^{2}$, then
\begin{eqnarray*}
&&E(Y_{t}^{m,n})^{2}+E\int_{t}^{T}\left| Z_{s}^{m,n}\right| ^{2}ds \\
&\leq &E[\xi ^{2}]+E\int_{t}^{T}\left| Y_{s}^{m,n}\right|
^{2}ds+E\int_{t}^{T}\left| f(s,0)\right| ^{2}ds+\frac{1}{\alpha }%
E[\sup_{0\leq t\leq T}(L_{t}^{+})^{2}]+\frac{1}{\alpha }E[\sup_{0\leq t\leq
T}(U_{t}^{-})^{2}] \\
&&+\alpha E[m\int_{t}^{T}(L_{s}-Y_{s}^{m,n})^{+}ds]^{2}+\alpha
E[n\int_{t}^{T}(U_{s}-Y_{s}^{m,n})^{-}ds]^{2},
\end{eqnarray*}
for some $\alpha >0$, in view of
\[
\int_{t}^{T}Y_{s}^{m,n}(L_{s}-Y_{s}^{m,n})^{+}ds=%
\int_{t}^{T}L_{s}(L_{s}-Y_{s}^{m,n})^{+}ds-%
\int_{t}^{T}((L_{s}-Y_{s}^{m,n})^{+})^{2}ds\leq
\int_{t}^{T}L_{s}(L_{s}-Y_{s}^{m,n})^{+}ds,
\]
and $\int_{t}^{T}Y_{s}^{m,n}(U_{s}-Y_{s}^{m,n})^{-}ds\leq
\int_{t}^{T}U_{s}(U_{s}-Y_{s}^{m,n})^{-}ds$. So
\begin{equation}
E\int_{t}^{T}\left| Z_{s}^{m,n}\right| ^{2}ds\leq C+\alpha
(E[m\int_{t}^{T}(L_{s}-Y_{s}^{m,n})^{+}ds]^{2}+E[n%
\int_{t}^{T}(U_{s}-Y_{s}^{m,n})^{-}ds)^{2}].  \label{est-zmn1}
\end{equation}

We need to prove that there exists a constant $C$ independent of $m,n$ such
that for any $0\leq t\leq T$%
\[
E[m\int_{t}^{T}(L_{s}-Y_{s}^{m,n})^{+}ds]^{2}+E[n%
\int_{t}^{T}(U_{s}-Y_{s}^{m,n})^{-}ds]^{2}\leq C+8E\int_{t}^{T}\left|
Z_{s}^{m,n}\right| ^{2}ds.
\]
In fact, let us consider the stopping time
\begin{eqnarray*}
\tau _{1} &=&\inf \{r\geq t|Y_{r}^{m,n}\geq U_{r}\}\wedge T,\sigma _{1}=\inf
\{r\geq \tau _{1}|Y_{r}^{m,n}=L_{r}\}\wedge T, \\
\tau _{2} &=&\inf \{r\geq \sigma _{1}|Y_{r}^{m,n}=U_{r}\}\wedge T,
\end{eqnarray*}
and so on. Since $L<U$ on $[0,T)$, and $L$ and $U$ are continuous, then when
$k\rightarrow \infty $, we have $\tau _{k}\nearrow T$, $\sigma _{k}\nearrow
T $. Obviously $Y^{m,n}\geq L$ on the interval $[\tau _{k},\sigma _{k}]$, so
we get
\[
Y_{\tau _{k}}^{m,n}=Y_{\sigma _{k}}^{m,n}+\int_{\tau _{k}}^{\sigma
_{k}}f(s,Y_{s}^{m,n})ds-n\int_{\tau _{k}}^{\sigma
_{k}}(Y_{s}^{m,n}-U_{s})^{+}ds-\int_{\tau _{k}}^{\sigma
_{k}}Z_{s}^{m,n}dB_{s}.
\]
On the other hand
\begin{eqnarray*}
Y_{\tau _{k}}^{m,n} &\geq &J_{\tau _{k}}\mbox{, on }\{\tau
_{k}<T\},Y_{\tau
_{k}}^{m,n}=J_{\tau _{k}}=\xi \mbox{, on }\{\tau _{k}=T\}, \\
Y_{\sigma _{k}}^{m,n} &\leq &J_{\sigma _{k}}\mbox{, on }\{\sigma
_{k}<T\},Y_{\sigma _{k}}^{m,n}=J_{\sigma _{k}}=\xi \mbox{, on
}\{\sigma _{k}=T\},
\end{eqnarray*}
and these inequalities imply that for all $k$, the following holds
\begin{eqnarray*}
n\int_{\tau _{k}}^{\sigma _{k}}(Y_{s}^{m,n}-U_{s})^{+}ds &\leq &J_{\sigma
_{k}}-J_{\tau _{k}}+\int_{\tau _{k}}^{\sigma
_{k}}f(s,Y_{s}^{m,n})ds-\int_{\tau _{k}}^{\sigma _{k}}Z_{s}^{m,n}dB_{s} \\
&\leq &\int_{\tau _{k}}^{\sigma _{k}}(\phi _{s}-Z_{s}^{m,n})dB_{s}+V_{\sigma
_{k}}^{-}-V_{\tau _{k}}^{-}+\int_{\tau _{k}}^{\sigma _{k}}\left|
f(s,Y_{s}^{m,n})\right| ds.
\end{eqnarray*}
Notice that on the interval $[\sigma _{k},\tau _{k+1}]$, $Y_{s}^{m,n}\leq
U_{s}$; we obtain by summing in $k$%
\[
n\int_{t}^{T}(Y_{s}^{m,n}-U_{s})^{+}ds\leq \int_{t}^{T}((\phi
_{s}-Z_{s}^{m,n})(\sum_{k}1_{[\tau _{k},\sigma
_{k})}(s))dB_{s}+V_{T}^{-}+\int_{t}^{T}\left| f(s,Y_{s}^{m,n})\right| ds.
\]

By squaring and taking the expectation, with (\ref{est-fmn}), we get
\begin{eqnarray}
&&E[n\int_{t}^{T}(Y_{s}^{m,n}-U_{s})^{+}ds]^{2}  \label{est-kmn1} \\
&\leq &4E\int_{t}^{T}\left| \phi _{s}\right| ^{2}ds+4E\int_{t}^{T}\left|
Z_{s}^{m,n}\right| ^{2}ds+2E[(V_{T}^{-})^{2}]+2E(\int_{t}^{T}\left|
f(s,Y_{s}^{m,n})\right| ds)^{2}  \nonumber \\
&\leq &C+4E\int_{t}^{T}\left| Z_{s}^{m,n}\right| ^{2}ds,  \nonumber
\end{eqnarray}
in the same way, we obtain
\begin{equation}
E[m\int_{t}^{T}(L_{s}-Y_{s}^{m,n})^{+}ds]^{2}\leq C+4E\int_{t}^{T}\left|
Z_{s}^{m,n}\right| ^{2}ds.  \label{est-kmn2}
\end{equation}

By (\ref{est-kmn1}) and (\ref{est-kmn2}), and (\ref{est-zmn1}), with $\alpha
=\frac{1}{16}$, it follows
\begin{equation}
E\int_{t}^{T}\left| Z_{s}^{m,n}\right| ^{2}ds\leq C,  \label{est-zmn2}
\end{equation}
then
\begin{equation}
E[(K_{T}^{m,n,+})^{2}+(K_{T}^{m,n,-})^{2}]\leq C.  \label{est-kmn3}
\end{equation}
$\square $

Let $m\rightarrow \infty $, due to the convergence results in step 1 of the
proof in \cite{LMX}, $Y^{m,n}\rightarrow Y^{n}$ in $\mathbf{S}^{2}(0,T)$, $%
K^{m,n,+}\rightarrow K^{n,+}$ in $\mathbf{A}^{2}(0,T)$, and $%
Z^{m,n}\rightarrow Z^{n}$ in $\mathbf{H}_{d}^{2}(0,T)$, where $%
(Y^{n},Z^{n},K^{n,+})$ is the solution of the one lower barrier RBSDE$(\xi
,f_{n},L)$, with $f_{n}(s,y)=f(s,y)-n(y-U_{s})^{+}$. So
\[
Y_{t}^{n}=\xi
+\int_{t}^{T}f(s,Y_{s}^{n})ds+K_{T}^{n,+}-K_{t}^{n,+}-n%
\int_{t}^{T}(Y_{s}^{n}-U_{s})^{+}ds-\int_{t}^{T}Z_{s}^{n}dB_{s},
\]
$Y_{t}^{n}\geq L_{t}$, $0\leq t\leq T$, $%
\int_{0}^{T}(Y_{s}^{n}-L_{s})dK_{s}^{n}=0$. Thank to the uniform
estimations, which we got as above, we know that there exists a constant $C$
independent of $n$ and $t$, s.t.
\begin{equation}
\sup_{0\leq t\leq T}(Y_{t}^{n})^{2}+f(t,Y_{t}^{n})\leq C,  \label{est-ynu}
\end{equation}
and
\begin{equation}
E\int_{0}^{T}\left| Z_{s}^{n}\right|
^{2}ds+E[(K_{T}^{n,+})^{2}]+E[(K_{T}^{n,-})^{2}]\leq C  \label{est-zkn}
\end{equation}
where $K_{t}^{n,-}=n\int_{0}^{t}(Y_{s}^{n}-U_{s})^{+}ds$. Then by the
comparison theorem 4.3 in \cite{LMX}, we deduce that $Y_{t}^{n}\searrow
Y_{t} $, for $t\in [0,T]$, as $n\rightarrow \infty $, and by the dominated
convergence theorem
\begin{equation}
E\int_{0}^{T}(Y_{s}^{n}-Y_{s})^{2}ds\rightarrow 0,\mbox{ as
}n\rightarrow \infty .  \label{con-ymn1}
\end{equation}

Then we want to prove the convergence of $(Z^{n})$ in $\mathbf{H}%
_{d}^{2}(0,T)$. For this, we need the following lemma, which is analogue as
Lemma 4 in \cite{LS}. With (\ref{est-fmn}), (\ref{est-mn}) and (\ref
{assup-bon}), we can easily get it, so we omit the proof.

\begin{lemma}
\label{con-uy}
\begin{equation}
\lim_{n\rightarrow \infty }E(\sup_{0\leq t\leq
T}((Y_{t}^{n}-U_{t})^{+})^{2}=0.  \label{con-uyk}
\end{equation}
\end{lemma}

For $n,p\in \mathbf{N}$, applying It\^{o}'s formula to $\left|
Y^{n}-Y^{p}\right| ^{2}$, and taking the expectation, then
\begin{eqnarray*}
&&E(Y_{t}^{n}-Y_{t}^{p})^{2}+E\int_{t}^{T}\left| Z_{s}^{n}-Z_{s}^{p}\right|
^{2}ds \\
&\leq
&2E\int_{t}^{T}(Y_{s}^{n}-U_{s})^{+}dK_{s}^{p,-}+2E%
\int_{t}^{T}(Y_{s}^{p}-U_{s})^{+}dK_{s}^{n,-} \\
&\leq &2(E[(\sup_{0\leq t\leq T}(Y_{s}^{n}-U_{s})^{+})^{2}])^{\frac{1}{2}%
}(E(K_{T}^{p,-})^{2})^{\frac{1}{2}}+2(E[(\sup_{0\leq t\leq
T}(Y_{s}^{p}-U_{s})^{+})^{2}])^{\frac{1}{2}}(E(K_{T}^{n,-})^{2})^{\frac{1}{2}%
},
\end{eqnarray*}
since $\int_{t}^{T}(Y_{s}^{n}-Y_{s}^{p})d(K_{s}^{n,+}-K_{s}^{p,+})\leq 0$.
So by (\ref{con-uyk}) and (\ref{est-zkn}), as $n,p\rightarrow \infty $, $%
E\int_{t}^{T}\left| Z_{s}^{n}-Z_{s}^{p}\right| ^{2}ds\rightarrow 0$, which
implies $\{Z^{n}\}$ is a Cauchy sequence in $\mathbf{H}_{d}^{2}(0,T)$. So
there exists a process $Z\in \mathbf{H}_{d}^{2}(0,T)$, s.t., as $%
n\rightarrow \infty $,
\[
E\int_{t}^{T}\left| Z_{s}^{n}-Z_{s}\right| ^{2}ds\rightarrow 0.
\]

Moreover by It\^{o}'s formula, we have
\begin{eqnarray*}
E[\sup_{0\leq t\leq T}\left| Y_{t}^{n}-Y_{t}^{p}\right| ^{2}] &\leq
&2E\int_{t}^{T}(Y_{s}^{n}-U_{s})^{+}dK_{s}^{p,-}+2E%
\int_{t}^{T}(Y_{s}^{p}-U_{s})^{+}dK_{s}^{n,-} \\
&&+2E[\sup_{0\leq t\leq T}\int_{t}^{T}\left| Z_{s}^{n}-Z_{s}^{p}\right|
\left| Y_{s}^{n}-Y_{s}^{p}\right| dB_{s}].
\end{eqnarray*}
By Burkholder-Davis-Gundy inequality and (\ref{con-uyk}), we get, as $%
n,p\rightarrow \infty $%
\[
E[\sup_{0\leq t\leq T}\left| Y_{t}^{n}-Y_{t}^{p}\right| ^{2}]\rightarrow 0,
\]
i.e. $Y^{n}\searrow Y$, in $\mathbf{S}^{2}(0,T).$

By the convergence of $Y_{t}^{n}$, i.e. $Y_{t}^{n}\searrow Y_{t}$, $0\leq
t\leq T$, and the fact that $f(s,y)$ is continuous and decreasing in $y$, we
get $f(s,Y_{s}^{n})\nearrow f(s,Y_{s})$, $0\leq s\leq T$. Moreover $\left|
f(s,Y_{s}^{n})\right| \leq C$. Using the monotonic convergence theorem, we
deduce that
\begin{equation}
E\int_{0}^{T}[f(t,Y_{t}^{n})-f(t,Y_{t})]^{2}dt\rightarrow 0,  \label{conv-f}
\end{equation}
i.e. the sequence $\{f(\cdot ,Y_{\cdot }^{n})\}$ is also a Cauchy sequence
in $\mathbf{H}^{2}(0,T)$.

Now we consider the convergence of the increasing processes $(K^{n,+})$ and $%
(K^{n,-})$. By the comparison theorem 4.3 in \cite{LMX}, we get $%
K_{t}^{n,+}\geq K_{t}^{p,+}$, $K_{t}^{n,+}-K_{s}^{n,+}\geq
K_{t}^{p,+}-K_{s}^{p,+}$, for $0\leq s\leq t\leq T$. So for $0\leq t\leq T$,
$K_{t}^{n,+}\nearrow K_{t}^{+}$, with $E[(K_{t}^{n,+})^{2}]\leq C$, we get
that $E[(K_{t}^{+})^{2}]\leq C$. Furthermore, $K_{T}^{n,+}-K_{T}^{p,+}\geq
K_{t}^{n,+}-K_{t}^{p,+}$, which follows
\[
E[\sup_{0\leq t\leq T}(K_{t}^{n,+}-K_{t}^{p,+})^{2}]\leq
E[(K_{T}^{n,+}-K_{T}^{p,+})^{2}]\rightarrow 0,
\]
so $K^{n,+}\rightarrow K^{+}$ in $\mathbf{A}^{2}(0,T)$. On the other hand,
since $(Y^{n},Z^{n},K^{n+},K^{n-})$ satisfies
\[
Y_{t}^{n}=\xi
+\int_{t}^{T}f(s,Y_{s}^{n})ds+K_{T}^{n,+}-K_{t}^{n,+}-(K_{T}^{n,-}-K_{t}^{n,-})-\int_{t}^{T}Z_{s}^{n}dB_{s},
\]
and we can rewrite it in the following form
\[
K_{t}^{n,-}=Y_{t}^{n}-Y_{0}^{n}+\int_{0}^{t}f(s,Y_{s}^{n})ds+K_{t}^{n,+}-%
\int_{0}^{t}Z_{s}^{n}dB_{s}.
\]
Without losing the generality, for $p<n$, with BDG inequality, we get
\begin{eqnarray*}
&&E[\sup_{0\leq t\leq T}(K_{t}^{n,-}-K_{t}^{p,-})^{2}] \\
&\leq &5E[\sup_{0\leq t\leq
T}(Y_{t}^{n}-Y_{t}^{p})^{2}]+5(Y_{0}^{n}-Y_{0}^{p})^{2}+5TE(%
\int_{0}^{T}(f(s,Y_{s}^{n})-f(s,Y_{s}^{p}))^{2}ds) \\
&&+5E[(K_{T}^{n,+}-K_{T}^{p,+})^{2}]+CE%
\int_{0}^{t}(Z_{s}^{n}-Z_{s}^{p})^{2}ds \\
&\rightarrow &0,
\end{eqnarray*}
i.e. there exists a process $K^{-}\in \mathbf{A}^{2}(0,T)$, s.t. $%
K^{n,-}\rightarrow K^{-}$ in $\mathbf{A}^{2}(0,T)$, and the limit $%
(Y,Z,K^{+},K^{-})$ satisfies
\[
Y_{t}=\xi
+\int_{t}^{T}f(s,Y_{s})ds+K_{T}^{+}-K_{t}^{+}-(K_{T}^{-}-K_{t}^{-})-%
\int_{t}^{T}Z_{s}dB_{s}.
\]

Since for $n\in \mathbf{N}$, $Y_{t}^{n}\geq L_{t}$, $0\leq t\leq T$, so $%
Y_{t}\geq L_{t}$. The last is to check (4) of definition \ref{Def-tb}. Since
$(Y^{n},K^{n,+},K^{n,-})$ tends to $(Y,K^{+},K^{-})$ uniformly in $t$ in
probability, then the measure $dK^{n,+}$ converges to $dK^{+}$ weakly in
probability, so
\[
\int_{0}^{T}(Y_{t}^{n}-L_{t})dK_{t}^{n,+}\rightarrow
\int_{0}^{T}(Y_{t}-L_{t})dK_{t}^{+},
\]
in probability as $n\rightarrow \infty $. Obviously $%
\int_{0}^{T}(Y_{t}-L_{t})dK_{t}^{+}\geq 0$, On the other hand, for each $%
n\in \mathbf{N}$, $\int_{0}^{T}(Y_{t}^{n}-L_{t})dK_{t}^{n,+}=0$. Hence
\[
\int_{0}^{T}(Y_{t}-L_{t})dK_{t}^{+}=0,\mbox{ a.s.}
\]
Similarly, we have $\int_{0}^{T}(Y_{t}-U_{t})dK_{t}^{-}=0$. Consequently the
triple $(Y,Z,K^{+},K^{-})$ is solution of the RBSDE$(\xi ,f,L,U)$, under the
assumptions (\ref{assup-bon}). $\square $

\paragraph{\textbf{Step 2. }}

In this step, we consider the case of a barrier $L$ which satisfies the
assumption 2.3'\textbf{-(i)}:
\[
E[\varphi ^{2}(\sup_{0\leq t\leq T}(L_{t})^{+})]<\infty ,
\]
and $L^{+}\in \mathbf{S}^{2}(0,T)$, but we still assume that for some $C>0$,
\begin{equation}
\left| \xi \right| +\sup_{0\leq t\leq T}\left| f(t,0)\right| +\sup_{0\leq
t\leq T}(U_{t})^{-}\leq C.  \label{assup-bon2}
\end{equation}

For $n\in \mathbf{N}$, set $L^{n}=L\wedge n$, then $\sup_{0\leq t\leq
T}(L_{t}^{n})^{+}\leq n$ and $L_{t}^{n}\leq L_{t}$; so assumption 2.3'-(ii),
(iii) are satisfied and by the step 1, we know that there exists a triple $%
(Y^{n},Z^{n},K^{n})$, with $K^{n}=K^{n,+}-K^{n,-}$, which satisfies
\begin{equation}
Y_{t}^{n}=\xi
+\int_{t}^{T}f(s,Y_{s}^{n})ds+K_{T}^{n,+}-K_{t}^{n,+}-(K_{T}^{n,-}-K_{t}^{n,-})-\int_{t}^{T}Z_{s}^{n}dB_{s},
\label{equ-2l}
\end{equation}
$L_{t}^{n}\leq Y_{t}^{n}\leq U_{t}$, $0\leq t\leq T$, and $%
\int_{0}^{T}(Y_{t}^{n}-L_{t}^{n})dK_{t}^{n,+}=%
\int_{0}^{T}(Y_{t}^{n}-U_{t})dK_{t}^{n,-}=0$.

Consider the solution $(\overline{Y},\overline{Z},\overline{K})$ of one
lower barrier RBSDE$(\xi ,f,L)$ and the solution $(\underline{Y},\underline{Z%
},\underline{K})$ of the super barrier RBSDE$(\xi ,f,U)$, in fact these two
equations can be considered as the following two barriers RBSDE$(\xi ,f,L,%
\overline{U})$ and RBSDE$(\xi ,f,\underline{L},U)$, where $\underline{L}%
=-\infty $, $\overline{U}=+\infty $. By the comparison theorem \ref{comprg},
it follows that $\underline{Y}_{t}\leq Y_{t}^{n}\leq \overline{Y}_{t}$, $%
0\leq t\leq T$. So
\begin{equation}
E[\sup_{0\leq t\leq T}\left| Y_{t}^{n}\right| ^{2}]\leq \max \{E[\sup_{0\leq
t\leq T}\left| \overline{Y}_{t}\right| ^{2},E[\sup_{0\leq t\leq T}\left|
\underline{Y}_{t}\right| ^{2}\}\leq C.  \label{est-2yln}
\end{equation}
Since $L_{t}^{n}\leq L_{t}^{n+1}$, $0\leq t\leq T$, thanks to the comparison
theorem \ref{comprg}, $Y_{t}^{n}\nearrow Y_{t}$, $0\leq t\leq T$. From the
above estimate and Fatou's lemma, we get
\begin{equation}
E[\sup_{0\leq t\leq T}(Y_{t})^{2}]\leq C.  \label{est-ygbr}
\end{equation}
And
\begin{equation}
E\int_{0}^{T}\left| Y_{t}^{n}-Y_{t}\right| ^{2}dt\rightarrow 0,\mbox{as }%
n\rightarrow \infty ,  \label{con-s2yb}
\end{equation}
follows from the dominated convergence theorem.

Notice that $f$ is decreasing on $y$, then $f(t,\overline{Y}_{t})\leq
f(t,Y_{t}^{n})\leq f(t,\underline{Y}_{t})$, $0\leq t\leq T$, and with the
integral property of $\underline{Y}$ and $\overline{Y}$, we have
\begin{equation}
E[(\int_{0}^{t}f(s,Y_{s}^{n})ds)^{2}]\leq \max \{E[(\int_{0}^{t}f(s,%
\overline{Y}_{s})ds)^{2}],E[(\int_{0}^{t}f(s,\underline{Y}%
_{s})ds)^{2}]\}\leq C.  \label{est-2lfn}
\end{equation}

In order to prove the convergence of $(Z^{n},K^{n})$, we first need a-priori
estimations. We apply the It\^{o} formula to $\left| Y_{t}^{n}\right| ^{2}$
on the interval $[t,T]$,
\begin{eqnarray}
&&\ \ E\left| Y_{t}^{n}\right| ^{2}+E\int_{t}^{T}\left| Z_{s}^{n}\right|
^{2}ds  \label{est-itogb} \\
\ \ &\leq &E\left| \xi \right| ^{2}+E\int_{t}^{T}\left| Y_{s}^{n}\right|
^{2}ds+E\int_{t}^{T}\left| f(s,0)\right| ^{2}ds+(\alpha +\beta
)E[\sup_{0\leq t\leq T}\left| Y_{t}^{n}\right| ^{2}]  \nonumber \\
&&+\frac{1}{\alpha }E[(K_{T}^{n,+}-K_{t}^{n,+})^{2}]+\frac{1}{\beta }%
E[(K_{T}^{n,-}-K_{t}^{n,-})^{2}],  \nonumber
\end{eqnarray}
where $K^{n}=K^{n,+}-K^{n,-}$. We first use the comparison theorem to
estimate $K^{n,-}$. Consider the linear RBSDE$(\xi ,f(s,L_{s}{}^{-}),L,U)$,
by existence results of \cite{CK}, we know there exists $(\widetilde{Y},%
\widetilde{Z},\widetilde{K}^{+},\widetilde{K}^{-})\in \mathbf{S}%
^{2}(0,T)\times \mathbf{H}_{d}^{2}(0,T)\times \mathbf{A}^{2}(0,T)\times
\mathbf{A}^{2}(0,T)$ satisfying
\begin{eqnarray*}
\widetilde{Y}_{t} &=&\xi +\int_{t}^{T}f(s,(L_{s})^{-})ds+\widetilde{K}%
_{T}^{+}-\widetilde{K}_{t}^{+}-(\widetilde{K}_{T}^{-}-\widetilde{K}%
_{t}^{-})-\int_{t}^{T}\widetilde{Z}_{s}dB_{s}, \\
L_{t} &\leq &\widetilde{Y}_{t}\leq U_{t},\int_{0}^{T}(\widetilde{Y}%
_{t}-L_{t})d\widetilde{K}_{t}^{+}=\int_{0}^{T}(\widetilde{Y}_{t}-U_{t})d%
\widetilde{K}_{t}^{-}=0.
\end{eqnarray*}
Then we have the following lemma, which will be proved in Appendix.

\begin{lemma}
\label{com-l}For $0\leq s\leq t\leq T$, $K_{t}^{n,-}-K_{s}^{n,-}\leq
\widetilde{K}_{t}^{-}-\widetilde{K}_{s}^{-}$, and $K_{T}^{n,-}\leq
\widetilde{K}_{T}^{-}$.
\end{lemma}

Now we have
\[
E[(K_{T}^{n,-})^{2}]\leq E[(\widetilde{K}_{T}^{-})^{2}]\leq C.
\]
We rewrite the RBSDE$(\xi ,f,L^{n},U)$ (\ref{equ-2l}),
\[
K_{T}^{n,+}-K_{t}^{n,+}=Y_{t}^{n}-\xi
-\int_{t}^{T}f(s,Y_{s}^{n})ds+(K_{T}^{n,-}-K_{t}^{n,-})+%
\int_{t}^{T}Z_{s}^{n}dB_{s},
\]
hence
\begin{eqnarray}
E(K_{T}^{n,+}-K_{t}^{n,+})^{2} &\leq &5E\left| Y_{t}^{n}\right|
^{2}+5E\left| \xi \right| ^{2}+5E(\int_{t}^{T}f(s,Y_{s}^{n})ds)^{2}
\label{est-2lkp} \\
&&+5E[(K_{T}^{n,-}-K_{t}^{n,-})^{2}]+5E\int_{t}^{T}\left| Z_{s}^{n}\right|
^{2}ds  \nonumber \\
&\leq &C+5E\int_{t}^{T}\left| Z_{s}^{n}\right| ^{2}ds.  \nonumber
\end{eqnarray}
Then we substitute (\ref{est-2lkp}) into (\ref{est-itogb}), set $\alpha =10$%
, $\beta =1$, and with (\ref{assup-bon2}) and (\ref{est-2yln}), it follows
\begin{equation}
E(K_{T}^{n,+})^{2}+E\int_{0}^{T}\left| Z_{s}^{n}\right| ^{2}ds\leq C.
\label{est-zgb}
\end{equation}

Now for $n,p\in \mathbf{N}$, $n\geq p$, then $L_{t}^{n}\geq L_{t}^{p}$, $%
0\leq t\leq T$. We apply the It\^{o}'s formula to $(\left|
Y_{t}^{n}-Y_{t}^{p}\right| ^{2})$ on the interval $[t,T]$, and take
expectation
\begin{eqnarray*}
E[\left| Y_{t}^{n}-Y_{t}^{p}\right| ^{2}]+E\int_{t}^{T}\left|
Z_{s}^{n}-Z_{s}^{p}\right| ^{2}ds &\leq
&2E\int_{t}^{T}(L_{s}^{n}-L_{s}^{p})dK_{s}^{n,+}-2E%
\int_{t}^{T}(L_{s}^{n}-L_{s}^{p})dK_{s}^{p,+} \\
&\leq &2E\int_{t}^{T}(L_{s}^{n}-L_{s}^{p})dK_{s}^{n,+},
\end{eqnarray*}
in view of $\int_{t}^{T}(Y_{s}^{n}-Y_{s}^{p})d(K_{s}^{n,-}-K_{s}^{p,-})\geq
0 $. Since $L_{t}-L_{t}^{n}\downarrow 0$, for each $t\in [0,T]$, and $%
L_{t}-L_{t}^{n}$ is continuous, by the Dini's theorem, the convergence holds
uniformly on the interval $[0,T]$, i.e.
\begin{equation}
E[\sup_{0\leq t\leq T}(L_{t}-L_{t}^{n})^{2}]\rightarrow 0\mbox{, as }%
n\rightarrow \infty .  \label{con-l}
\end{equation}
Then with (\ref{est-2lkp}),
\begin{eqnarray*}
E\int_{0}^{T}\left| Z_{s}^{n}-Z_{s}^{p}\right| ^{2}ds &\leq &2(E(\sup_{0\leq
t\leq T}(L_{s}^{n}-L_{s}^{p})^{2})^{\frac{1}{2}}(E[(K_{T}^{n,+})^{2}])^{%
\frac{1}{2}} \\
&\leq &C(E(\sup_{0\leq t\leq T}(L_{s}^{n}-L_{s}^{p})^{2}])^{\frac{1}{2}%
}\rightarrow 0,
\end{eqnarray*}
as $n,p\rightarrow \infty $, i.e. $\{Z^{n}\}$ is a Cauchy sequence in the
space $\mathbf{H}_{d}^{2}(0,T)$, and there exists a process $Z\in \mathbf{H}%
_{d}^{2}(0,T)$, s.t. as $n\rightarrow \infty $,
\begin{equation}
E\int_{0}^{T}\left| Z_{s}^{n}-Z_{s}\right| ^{2}ds\rightarrow 0.
\label{con-znb}
\end{equation}

Furthermore from It\^{o}'s formula, we have
\begin{eqnarray*}
&&\ \ \sup_{0\leq t\leq T}\left| Y_{t}^{n}-Y_{t}^{p}\right| ^{2} \\
\  &\leq &2\sup_{0\leq t\leq
T}\int_{t}^{T}(L_{s}^{n}-L_{s}^{p})d(K_{s}^{n,+}-K_{s}^{p,+})+2\sup_{0\leq
t\leq T}\left|
\int_{t}^{T}(Y_{s}^{n}-Y_{s}^{p})(Z_{s}^{n}-Z_{s}^{p})dB_{s}\right| .
\end{eqnarray*}
Taking the expectation on the both sides, by BDG inequality and (\ref
{est-zgb}), we get
\begin{eqnarray*}
&&\ \ E\sup_{0\leq t\leq T}\left| Y_{t}^{n}-Y_{t}^{p}\right| ^{2} \\
\  &\leq
&2E[\int_{0}^{T}(L_{s}^{n}-L_{s}^{p})dK_{s}^{n,+}]+CE%
\int_{0}^{T}(Y_{s}^{n}-Y_{s}^{p})^{2}(Z_{s}^{n}-Z_{s}^{p})^{2}ds \\
\ \  &\leq &C(E[\sup_{0\leq t\leq T}(L_{s}^{n}-L_{s}^{p})^{2}])^{\frac{1}{2}%
}+\frac{1}{2}E\sup_{0\leq t\leq T}\left| Y_{s}^{n}-Y_{s}^{p}\right|
^{2}+CE\int_{0}^{T}\left| Z_{s}^{n}-Z_{s}^{p}\right| ^{2}ds.
\end{eqnarray*}
Hence, by (\ref{con-znb}) and (\ref{con-l}), as $n,p\rightarrow \infty $,
\begin{equation}
E\sup_{0\leq t\leq T}\left| Y_{t}^{n}-Y_{t}^{p}\right| ^{2}\rightarrow 0,
\label{cauchy-ynb}
\end{equation}
i.e. $\{Y^{n}\}$ is a Cauchy sequence in the space $\mathbf{S}^{2}(0,T)$,
which implies that there exists a process $Y\in \mathbf{S}^{2}(0,T)$, s.t.
as $n\rightarrow \infty $,
\begin{equation}
E\sup_{0\leq t\leq T}\left| Y_{t}^{n}-Y_{t}\right| ^{2}\rightarrow 0.
\label{con-ynb}
\end{equation}
Moreover, since $f$ is continuous and decreasing on $y$, with $%
Y_{t}^{n}\nearrow Y_{t}$,
\[
f(t,Y_{t}^{n})-f(t,Y_{t})\searrow 0,\mbox{ }0\leq t\leq T\mathbf{.}
\]
By the monotonic limit theorem, we get $%
\int_{0}^{T}[f(t,Y_{t}^{n})-f(t,Y_{t})]dt\searrow 0$, and with (\ref
{est-2lfn}), it follows $E[(\int_{0}^{T}f(t,Y_{t})dt)^{2}]\leq C$, then
\begin{equation}
E[(\int_{0}^{T}(f_{n}(t,Y_{t}^{n})-f(t,Y_{t}))dt)^{2}]\rightarrow 0,
\label{con-2lfnb}
\end{equation}
as $n\rightarrow \infty $.

From corollary \ref{compr2}, we know that for $\forall t\in [0,T]$, $%
K_{t}^{n,-}$ is increasing with respect to $n$, and with $%
E[(K_{t}^{n,-})^{2}]\leq C$, there exists $K_{t}^{-}$ such that $%
K_{t}^{n,-}\nearrow K_{t}^{-}$ in $\mathbf{L}^{2}(\mathcal{F}_{t})$. Since
for each $t\in [0,T]$, $E[(K_{t}^{n,+})^{2}]\leq C$, the sequence $%
(K_{t}^{n,+})$ has weak limit $K_{t}^{+}$ in $\mathbf{L}^{2}(\mathcal{F}%
_{t}) $, with $E[(K_{t}^{+})^{2}]\leq C$. Then for $0\leq t\leq T$, $%
(Y,Z,K^{+},K^{-})$ satisfies
\begin{equation}
Y_{t}=\xi
+\int_{t}^{T}f(s,Y_{s})ds+K_{T}^{+}-K_{t}^{+}-(K_{T}^{-}-K_{t}^{-})-%
\int_{t}^{T}Z_{s}dB_{s}.  \label{solutionb}
\end{equation}

We will then prove that the convergence of $\{K^{n,+}\}$ and $\{K^{n,-}\}$
also holds in strong sense. First, we consider $\{K^{n,-}\}$, for $n$, $p\in
\mathbf{N}$, with $n\geq p$, since $L_{t}^{n}\geq L_{t}^{p}$, by corollary
\ref{compr2}, we have for $0\leq s\leq t\leq T$, $K_{t}^{n,-}-K_{s}^{n,-}%
\geq K_{t}^{p,-}-K_{s}^{p,-}$. So $0\leq K_{t}^{n,-}-K_{t}^{p,-}\leq
K_{T}^{n,-}-K_{T}^{p,-}$, and it follows immediately by letting $%
n\rightarrow \infty $%
\[
0\leq K_{t}^{-}-K_{t}^{p,-}\leq K_{T}^{-}-K_{T}^{p,-}.
\]
This inequality yields as $p\rightarrow \infty $,
\begin{equation}
E\sup_{0\leq t\leq T}\left| K_{t}^{-}-K_{t}^{p,-}\right| ^{2}\leq E\left|
K_{T}^{-}-K_{T}^{p,-}\right| ^{2}\rightarrow 0.  \label{con-knb}
\end{equation}
Then we consider the term $\{K^{n,+}\}$. For this we rewrite (\ref{equ-2l})
and (\ref{solutionb}) in the forward form:
\begin{eqnarray*}
K_{t}^{n,+}
&=&Y_{0}^{n}-Y_{t}^{n}-\int_{0}^{t}f(s,Y_{s}^{n})ds+K_{t}^{n,-}+%
\int_{0}^{t}Z_{s}^{n}dB_{s} \\
K_{t}^{+}
&=&Y_{0}-Y_{t}-\int_{0}^{t}f(s,Y_{s})ds+K_{t}^{-}+\int_{0}^{t}Z_{s}dB_{s},
\end{eqnarray*}
so consider the difference and take expectation on the both sides, by the
BDG inequality, and $f(s,Y_{s}^{n})\geq f(s,Y_{s})$, it follows
\begin{eqnarray*}
E[\sup_{0\leq t\leq T}\left| K_{t}^{n,+}-K_{t}^{+}\right| ^{2}] &\leq
&5\left| Y_{0}^{n}-Y_{0}\right| ^{2}+5E[\sup_{0\leq t\leq T}\left|
Y_{t}^{n}-Y_{t}\right|
^{2}]+5E(\int_{0}^{T}[f(s,Y_{s}^{n})-f(s,Y_{s})]ds)^{2} \\
&&\ +5E\sup_{0\leq t\leq T}\left| K_{t}^{n,-}-K_{t}^{-}\right|
^{2}+CE\int_{0}^{T}\left| Z_{s}^{n}-Z_{s}\right| ^{2}ds.
\end{eqnarray*}
Then by (\ref{con-ynb}), (\ref{con-2lfnb}), (\ref{con-knb}) and (\ref
{con-znb}), we deduce that
\[
E[\sup_{0\leq t\leq T}\left| K_{t}^{n,+}-K_{t}^{+}\right| ^{2}]\rightarrow
0.
\]

The last thing to check is that (3) and (4) are also satisfied. Since for
each $n\in \mathbf{N}$, $L_{t}^{n}\leq Y_{t}^{n}\leq U_{t}$, $0\leq t\leq T$%
, with $Y_{t}^{n}\nearrow Y_{t}$ and $L_{t}^{n}\nearrow L_{t}$, then $%
L_{t}\leq Y_{t}\leq U_{t}$. From another part, the processes $K^{n,+}$ and $%
K^{n,-}$ are increasing, so the limit $K^{+}$ and $K^{-}$ are also
increasing. Notice that $(Y^{n},K^{n,+},K^{n,-})$ tends to $(Y,K^{+},K^{-})$
uniformly in $t$ in probability, so the measure $dK^{n,+}$(resp. $dK^{n,-}$)
converges to $dK^{+}$ (resp. $dK^{-}$) weakly in probability. So
\[
\int_{0}^{T}(Y_{t}-L_{t})dK_{t}^{n,+}\rightarrow
\int_{0}^{T}(Y_{t}-L_{t})dK_{t}^{+}\mbox{ ,}%
\int_{0}^{T}(Y_{t}^{n}-U_{t})dK_{t}^{n,-}\rightarrow
\int_{0}^{T}(Y_{t}-U_{t})dK_{t}^{-},
\]
in probability as $n\rightarrow \infty $. Obviously $%
\int_{0}^{T}(Y_{t}-U_{t})dK_{t}^{-}\leq 0$. On the other hand, for each $%
n\in \mathbf{N}$, $\int_{0}^{T}(Y_{t}^{n}-U_{t})dK_{t}^{n,-}=0$. Hence
\[
\int_{0}^{T}(Y_{t}-U_{t})dK_{t}^{-}=0,\mbox{ a.s.}
\]
For the lower barrier, since $L^{n}$ converges to $L$ in $\mathbf{S}^{2}$,
as $n\rightarrow \infty $, we have
\begin{eqnarray*}
&&\
E\int_{0}^{T}(Y_{t}^{n}-L_{t}^{n})dK_{t}^{n,+}-E%
\int_{0}^{T}(Y_{t}-L_{t})dK_{t}^{+} \\
&=&E\int_{0}^{T}(Y_{t}^{n}-Y_{t})dK_{t}^{n,+}+E%
\int_{0}^{T}(Y_{t}-L_{t})d(K_{t}^{n,+}-K_{t}^{+})+E%
\int_{0}^{T}(L_{t}-L_{t}^{n})dK_{t}^{n,+} \\
\ &\leq &C(E[\sup_{0\leq t\leq T}(Y_{t}^{n}-Y_{t})^{2}])^{\frac{1}{2}%
}+E\int_{0}^{T}(Y_{t}-L_{t})d(K_{t}^{n,+}-K_{t}^{+})+C(E[\sup_{0\leq t\leq
T}(L_{t}-L_{t}^{n})^{2}])^{\frac{1}{2}} \\
\ &\rightarrow &0.
\end{eqnarray*}
Since $Y_{t}\geq L_{t}$, then $\int_{0}^{T}(Y_{t}-L_{t})dK_{t}^{+}\geq 0$,
while $E\int_{0}^{T}(Y_{t}^{n}-L_{t}^{n})dK_{t}^{n,+}=0$, so $%
E\int_{0}^{T}(Y_{t}-L_{t})dK_{t}^{+}=0$, then $%
\int_{0}^{T}(Y_{t}-L_{t})dK_{t}^{+}=0.\square $

\paragraph{\textbf{Step 3.}}

In this step, we study the general case for $L$ and $U$, when assumption
2.3' is satisfied:
\[
E[\varphi ^{2}(\sup_{0\leq t\leq T}(L_{t})^{+})]+E[\varphi ^{2}(\sup_{0\leq
t\leq T}(U_{t})^{-})]<\infty ,
\]
$L^{+}$, $U^{-}\in \mathbf{S}^{2}(0,T)$. But we still assume that for some $%
C>0$,
\begin{equation}
\left| \xi \right| +\sup_{0\leq t\leq T}\left| f(t,0)\right| \leq C.
\label{assup-bon3}
\end{equation}

For $n\in \mathbf{N}$, set $U^{n}=U\vee (-n)$; then $\sup_{0\leq t\leq
T}(U_{t}^{n})^{-}\leq n$ and $U^{n}\geq U$, so assumption 2.3'-(ii), (iii)
are satisfied, and by the step 2, we know that there exists a triple $%
(Y^{n},Z^{n},K^{n})$, with $K^{n}=K^{n,+}-K^{n,-}$, which satisfies
\begin{equation}
Y_{t}^{n}=\xi
+\int_{t}^{T}f(s,Y_{s}^{n})ds+K_{T}^{n,+}-K_{t}^{n,+}-(K_{T}^{n,-}-K_{t}^{n,-})-\int_{t}^{T}Z_{s}^{n}dB_{s},
\label{equ-ub}
\end{equation}
$L_{t}\leq Y_{t}^{n}\leq U_{t}^{n}$, $0\leq t\leq T$, and $%
\int_{0}^{T}(Y_{t}^{n}-L_{t})dK_{t}^{n,+}=%
\int_{0}^{T}(Y_{t}^{n}-U_{t}^{n})dK_{t}^{n,-}=0$.

Like in step 2, we consider the solution $(\overline{Y},\overline{Z},%
\overline{K})$ of the one lower barrier RBSDE$(\xi ,f,L)$ and the solution $(%
\underline{Y},\underline{Z},\underline{K})$ of the one super barrier RBSDE$%
(\xi ,f,U)$. Then by the comparison theorem \ref{comprg}, it follows that $%
\underline{Y}_{t}\leq Y_{t}^{n}\leq \overline{Y}_{t}$, $0\leq t\leq T$. So
\begin{equation}
E[\sup_{0\leq t\leq T}\left| Y_{t}^{n}\right| ^{2}]\leq \max \{E[\sup_{0\leq
t\leq T}\left| \overline{Y}_{t}\right| ^{2}],E[\sup_{0\leq t\leq T}\left|
\underline{Y}_{t}\right| ^{2}]\}\leq C.  \label{est-yupb}
\end{equation}
Since $U_{t}^{n+1}\leq U_{t}^{n}$, $0\leq t\leq T$, thanks to the comparison
theorem \ref{comprg}, $Y_{t}^{n}\searrow Y_{t}$, $0\leq t\leq T$. From (\ref
{est-yupb}) and Fatou's lemma, we get
\begin{equation}
E[\sup_{0\leq t\leq T}(Y_{t})^{2}]\leq C,  \label{est-yubu}
\end{equation}
and
\begin{equation}
E\int_{0}^{T}\left| Y_{t}^{n}-Y_{t}\right| ^{2}dt\rightarrow 0,\mbox{as }%
n\rightarrow \infty ,  \label{con-s2yub}
\end{equation}
which follows from the dominated convergence theorem.

Notice that $f$ is decreasing on $y$, then $f(t,\overline{Y}_{t})\leq
f(t,Y_{t}^{n})\leq f(t,\underline{Y}_{t})$, $0\leq t\leq T$, and with the
integral property of $\underline{Y}$ and $\overline{Y}$, we have
\begin{equation}
E[(\int_{0}^{t}f(s,Y_{s}^{n})ds)^{2}]\leq \max \{E[(\int_{0}^{t}f(s,%
\overline{Y}_{s})ds)^{2}],E[(\int_{0}^{t}f(s,\underline{Y}%
_{s})ds)^{2}]\}\leq C.  \label{est-fu}
\end{equation}
Then we use again the comparison theorem for the estimation of $K_{t}^{n,+}$%
. Consider the linear RBSDE$(\xi ,f(s,U_{s}{}^{+}),L,U)$, by results of \cite
{CK}, we know that there exists $(\widetilde{Y},\widetilde{Z},\widetilde{K}%
^{+},\widetilde{K}^{-})\in \mathbf{S}^{2}(0,T)\times \mathbf{H}%
_{d}^{2}(0,T)\times \mathbf{A}^{2}(0,T)\times \mathbf{A}^{2}(0,T)$
satisfying the following:
\begin{eqnarray*}
\widetilde{Y}_{t} &=&\xi +\int_{t}^{T}f(s,(U_{s})^{+})ds+\widetilde{K}%
_{T}^{+}-\widetilde{K}_{t}^{+}-(\widetilde{K}_{T}^{-}-\widetilde{K}%
_{t}^{-})-\int_{t}^{T}\widetilde{Z}_{s}dB_{s}, \\
L_{t} &\leq &\widetilde{Y}_{t}\leq U_{t},\int_{0}^{T}(\widetilde{Y}%
_{t}-L_{t})d\widetilde{K}_{t}^{+}=\int_{0}^{T}(\widetilde{Y}_{t}-U_{t})d%
\widetilde{K}_{t}^{-}=0.
\end{eqnarray*}
We admits for a instant the following lemma, which will be proved later.

\begin{lemma}
\label{com-bu}For $0\leq s\leq t\leq T$, $K_{t}^{n,+}-K_{s}^{n,+}\leq
\widetilde{K}_{t}^{+}-\widetilde{K}_{s}^{+}$, and $K_{T}^{n,+}\leq
\widetilde{K}_{T}^{+}$.
\end{lemma}

Now we have
\[
E[(K_{T}^{n,+})^{2}]\leq E[(\widetilde{K}_{T}^{+})^{2}]\leq C,
\]
Then apply the It\^{o}'s formula to $\left| Y_{t}^{n}\right| ^{2}$ on the
interval $[t,T]$, by the same method as in step 2, we have the following
estimates
\begin{equation}
E[(K_{T}^{n,-})^{2}]+E\int_{0}^{T}\left| Z_{s}^{n}\right| ^{2}ds\leq C.
\label{est-zkub}
\end{equation}

Since $U_{t}^{n}-U_{t}\downarrow 0$, for each $t\in [0,T]$, and $%
U_{t}^{n}-U_{t}$ is continuous, by the Dini's theorem again, the convergence
holds uniformly on the interval $[0,T]$, i.e.
\begin{equation}
E[\sup_{0\leq t\leq T}(U_{t}^{n}-U_{t})^{2}]\rightarrow 0\mbox{, as }%
n\rightarrow \infty .  \label{con-uu}
\end{equation}
Now we are in the same situation as step 2. With the same arguments, we
deduce that there exists processes $Y\in \mathbf{S}^{2}(0,T)$, $Z\in \mathbf{%
H}_{d}^{2}(0,T)$, $K^{+}\in \mathbf{A}^{2}(0,T)$, $K^{-}\in \mathbf{A}%
^{2}(0,T)$, s.t. as $n\rightarrow \infty $,
\[
E[\sup_{0\leq t\leq T}\left| Y_{t}^{n}-Y_{t}\right| +\int_{0}^{T}\left|
Z_{s}^{n}-Z_{s}\right| ^{2}ds+\sup_{0\leq t\leq T}\left|
K_{t}^{n,+}-K_{t}^{+}\right| ^{2}+\sup_{0\leq t\leq T}\left|
K_{t}^{n,-}-K_{t}^{-}\right| ^{2}]\rightarrow 0,
\]
which satisfies
\begin{equation}
Y_{t}=\xi
+\int_{t}^{T}f(s,Y_{s})ds+K_{T}^{+}-K_{t}^{+}-(K_{T}^{-}-K_{t}^{-})-%
\int_{t}^{T}Z_{s}dB_{s}.  \label{solutionub}
\end{equation}

The last thing to check is that (3) and (4) of definition \ref{Def-tb} are
satisfied. Since for each $n\in \mathbf{N}$, $L_{t}\leq Y_{t}^{n}\leq
U_{t}^{n}$, $0\leq t\leq T$, with $Y_{t}^{n}\searrow Y_{t}$ and $%
U_{t}^{n}\searrow U_{t}$, then $L_{t}\leq Y_{t}\leq U_{t}$. On the other
hand, the processes $K^{n,+}$ and $K^{n,-}$ are increasing, so the limit $%
K^{+}$ and $K^{-}$ are also increasing. Notice that $(Y^{n},K^{n,+},K^{n,-})$
tends to $(Y,K^{+},K^{-})$ uniformly in $t$ in probability, and $U^{n}$
converges to $U$ in $\mathbf{S}^{2}$, as $n\rightarrow \infty $, similarly
as step 2, we get
\[
\int_{0}^{T}(Y_{t}-K_{t})dK_{t}^{+}=\int_{0}^{T}(Y_{t}-U_{t})dK_{t}^{-}=0,%
\mbox{ a.s.}
\]
$\square $

\paragraph{\textbf{Step 4.}}

In this step, we will partly relax the bounded assumption for $\xi $ and $%
f(t,0)$. We only suppose that for a constant $c$,
\begin{equation}
\xi \geq c\mbox{ and }\inf_{0\leq t\leq T}f(t,0)\geq c.
\label{assup2-blow}
\end{equation}
We approximate $\xi $ and $f(t,0)$ by a sequence whose elements satisfy the
bounded assumption in step 3, as following: for $n\in \mathbf{N}$, set
\[
\xi _{n}=\xi \wedge n,f_{n}(t,y)=f(t,y)-f(t,0)+f(t,0)\wedge n.
\]
Obviously, $(\xi ^{n},f^{n})$ satisfies the assumptions of the step 3, and
since $\xi \in \mathbf{L}^{2}(\mathcal{F}_{T})$, $f(t,0)\in \mathbf{H}%
^{2}(0,T)$, then
\begin{equation}
E[\left| \xi ^{n}-\xi \right| ^{2}]\rightarrow 0,E\int_{0}^{T}\left|
f(t,0)-f_{n}(t,0)\right| ^{2}\rightarrow 0,  \label{con2-ini}
\end{equation}
as $n\rightarrow \infty $.

From the results in step 3, for each $n\in \mathbf{N}$, there exists $%
(Y_{t}^{n},Z_{t}^{n},K_{t}^{n})_{0\leq t\leq T}\in \mathbf{S}^{2}(0,T)\times
\mathbf{H}_{d}^{2}(0,T)\times \mathbf{VF}^{2}(0,T)$, with $%
K^{n}=K^{n,+}-K^{n,-}$, which is the unique solution of the RBSDE$(\xi
^{n},f_{n},L,U)$, i.e.
\begin{eqnarray}
Y_{t}^{n} &=&\xi
^{n}+%
\int_{t}^{T}f_{n}(s,Y_{s}^{n})ds+K_{T}^{n,+}-K_{t}^{n,+}-(K_{T}^{n,-}-K_{t}^{n,-})-\int_{t}^{T}Z_{s}^{n}dB_{s},
\label{equa} \\
L_{t} &\leq &Y_{t}^{n}\leq
U_{t},\int_{0}^{T}(Y_{t}^{n}-L_{t})dK_{t}^{n,+}=%
\int_{0}^{T}(Y_{t}^{n}-U_{t})dK_{t}^{n,-}=0.  \nonumber
\end{eqnarray}
Like in step 3, we consider the solution $(\overline{Y},\overline{Z},%
\overline{K})$ of one lower barrier RBSDE$(\xi ,f,L)$ and the solution $(%
\underline{Y},\underline{Z},\underline{K})$ of one super barrier RBSDE$(\xi
^{-},\underline{f},U)$, where $\xi ^{-}$ is the negative part of $\xi $, $%
\underline{f}(t,y)=f(t,y)-f(t,0)+(f(t,0))^{-}$. Then we can take the RBSDE$%
(\xi ,f,L)$ (resp. RBSDE$(\xi ^{-},\underline{f},U)$) as a RBSDE with two
barriers associated to the parameters $(\xi ,f,L,\overline{U})$ (resp. $(\xi
^{-},\underline{f},\underline{L},U)$), where $\overline{U}=\infty $ and $%
\underline{L}=-\infty $. By the comparison theorem \ref{comprg}, since
\[
\xi \geq \xi ^{n}\geq \xi ^{-},f(t,y)\geq f_{n}(t,y)\geq \underline{f}(t,y),
\]
it follows that
\[
\underline{Y}_{t}\leq Y_{t}^{n}\leq \overline{Y}_{t},0\leq t\leq T.
\]
So
\[
E[\sup_{0\leq t\leq T}\left| Y_{t}^{n}\right| ^{2}]\leq \max \{E[\sup_{0\leq
t\leq T}\left| \overline{Y}_{t}\right| ^{2},E[\sup_{0\leq t\leq T}\left|
\underline{Y}_{t}\right| ^{2}\}\leq C.
\]
Then by the comparison theorem \ref{compr4}, since for all $(s,y)\in
[0,T]\times \mathbb{R}$, $n\in \mathbf{N}$, $\xi _{1}\leq \xi _{n}$, $%
f_{1}(s,y)\leq f_{n}(s,y)$, we have $K_{t}^{1,+}\geq K_{t}^{n,+}\geq 0$ for $%
0\leq t\leq T$, so $E[(K_{t}^{n,+})^{2}]\leq E[(K_{t}^{1,+})^{2}]\leq C$.
Following the same steps, we deduce that
\begin{equation}
E[\int_{0}^{t}f(s,Y_{s}^{n})ds)^{2}]+E\int_{0}^{T}\left| Z_{s}^{n}\right|
^{2}ds+E[(K_{t}^{n,-})^{2}]+E[(K_{t}^{n,+})^{2}]\leq C.  \label{est2-bfkz}
\end{equation}
Due to the comparison theorem \ref{comprg}, since for all $(s,y)\in
[0,T]\times \mathbb{R}$, $n\in \mathbf{N}$, $\xi _{n}\leq \xi _{n+1}$, $%
f_{n}(s,y)\leq f_{n+1}(s,y)$, we have $Y_{t}^{n}\leq Y_{t}^{n+1}$, $0\leq
t\leq T$, a.s. Hence
\begin{equation}
Y_{t}^{n}\nearrow Y_{t},0\leq t\leq T.\mbox{ a.s.}  \label{con2-mon}
\end{equation}

Applying It\^{o} formula to $\left| Y_{t}^{n}-Y_{t}^{p}\right| ^{2}$, for $%
n,p\in \mathbf{N}$, $n\geq p$, on $[t,T]$, we get
\begin{eqnarray*}
&&E\left| Y_{t}^{n}-Y_{t}^{p}\right| ^{2}+E\int_{t}^{T}\left|
Z_{s}^{n}-Z_{s}^{p}\right| ^{2}ds \\
&\leq &E\left| \xi ^{n}-\xi ^{p}\right| ^{2}+E\int_{t}^{T}\left|
Y_{s}^{n}-Y_{s}^{p}\right| ^{2}ds+E\int_{t}^{T}\left|
f_{n}(s,0)-f_{p}(s,0)\right| ^{2}ds,
\end{eqnarray*}
since $\int_{t}^{T}(Y_{s}^{n}-Y_{s}^{p})d(K_{s}^{n,+}-K_{s}^{p,+})-%
\int_{t}^{T}(Y_{s}^{n}-Y_{s}^{p})d(K_{s}^{n,-}-K_{s}^{p,-})\leq 0$. Hence
from Gronwall's inequality and (\ref{con2-ini}), we deduce
\begin{equation}
\sup_{0\leq t\leq T}E\left| Y_{t}^{n}-Y_{t}^{p}\right| ^{2}\rightarrow
0,\;\;E\int_{0}^{T}\left| Z_{s}^{n}-Z_{s}^{p}\right| ^{2}ds\rightarrow 0.
\label{cauchy2-ynz}
\end{equation}
Consequently there exists $(Z_{t})_{0\leq t\leq T}\in \mathbf{H}%
_{d}^{2}(0,T) $, s.t.
\begin{equation}
E\int_{0}^{T}\left| Z_{s}^{n}-Z_{s}\right| ^{2}ds\rightarrow 0.
\label{con2-zn2}
\end{equation}

Using again It\^{o} formula, taking $\sup $ and the expectation, in view of
the BDG inequality, $Y_{t}^{n}\geq Y_{t}^{p}$, assumption 2.4-(iii) and $%
f_{n}(t,0)\geq f_{p}(t,0)$, we get
\begin{eqnarray*}
E[\sup_{0\leq t\leq T}\left| Y_{t}^{n}-Y_{t}^{p}\right| ^{2}] &\leq &E\left|
\xi ^{n}-\xi ^{p}\right| +4TE\int_{0}^{T}\left| f_{n}(s,0)-f_{p}(s,0)\right|
^{2}ds+\frac{1}{4}E\sup_{0\leq t\leq T}\left| Y_{s}^{n}-Y_{s}^{p}\right| ^{2}
\\
&&+\frac{1}{4}E[\sup_{0\leq t\leq T}\left| Y_{t}^{n}-Y_{t}^{p}\right|
^{2}]+cE\int_{0}^{T}\left| Z_{s}^{n}-Z_{s}^{p}\right| ^{2}ds.
\end{eqnarray*}
From (\ref{con2-ini}) and (\ref{cauchy2-ynz}), it follows $E[\sup_{0\leq
t\leq T}\left| Y_{t}^{n}-Y_{t}^{p}\right| ^{2}]\rightarrow 0$, as $%
n,p\rightarrow \infty $, i.e. the sequence $\{Y^{n}\}$ is a Cauchy sequence
in the space $\mathbf{S}^{2}(0,T)$. Consequently, with (\ref{con2-mon}), we
have $Y\in \mathbf{S}^{2}(0,T)$ and
\begin{equation}
E[\sup_{0\leq t\leq T}\left| Y_{t}^{n}-Y_{t}\right| ^{2}]\rightarrow 0.
\label{con2-yn2}
\end{equation}

By the comparison theorem \ref{compr4}, since for all $(s,y)\in
[0,T]\times \mathbb{R}$, $n\in \mathbf{N}$, $\xi _{n}\leq \xi
_{n+1}$, $f_{n}(s,y)\leq f_{n+1}(s,y)$, we have $K_{t}^{n,+}\geq
K_{t}^{n+1,+}\geq 0$, and $0\leq K_{t}^{n,-}\leq K_{t}^{n+1,-}$ for
$0\leq t\leq T$, so
\begin{equation}
K_{t}^{n,+}\searrow K_{t}^{+},\mbox{ }K_{t}^{n,-}\nearrow
K_{t}^{-},\mbox{ } \label{con-kn}
\end{equation}
with (\ref{est2-bfkz}), by the monotonic limit theorem, it follows that $%
K_{t}^{n,+}\rightarrow K_{t}^{+}$, $K_{t}^{n,-}\rightarrow K_{t}^{-}$ in $%
\mathbf{L}^{2}(\mathcal{F}_{t})$, and $E[(K_{t}^{+})^{2}+(K_{t}^{-})^{2}]<%
\infty $, moreover, $(K_{t}^{+})_{0\leq t\leq T}$ and $(K_{t}^{-})_{0\leq
t\leq T}$ are increasing.

Notice that since $f(t,y)$ is decreasing and continuous in $y$, and $%
Y_{t}^{n}\nearrow Y_{t}$, we have $f(t,Y_{t}^{n})\searrow f(t,Y_{t})$. Then
by the monotonic limit theorem, $\int_{0}^{t}f(s,Y_{s}^{n})ds\searrow
\int_{0}^{t}f(s,Y_{s})ds$. With (\ref{est2-bfkz}), it follows that $%
\int_{0}^{t}f(s,Y_{s}^{n})ds\rightarrow \int_{0}^{t}f(s,Y_{s})ds$ in $%
\mathbf{L}^{2}(\mathcal{F}_{t})$, as $n\rightarrow \infty $.

Now we need to prove that the convergence of $\{K^{n,+}\}$ and
$\{K^{n,-}\}$ holds in a stronger sense. Using again the comparison
theorem \ref{compr4}, since for all $(s,y)\in [0,T]\times
\mathbb{R}$, $n$, $p\in \mathbf{N}$, with $n\geq p$, $\xi _{p}\leq
\xi _{n}$, $f_{p}(s,y)\leq f_{n}(s,y)$, we have for $0\leq s\leq
t\leq T$,
\[
K_{t}^{p,+}-K_{s}^{p,+}\geq K_{t}^{n,+}-K_{s}^{n,+}\geq 0,
\]
Then let $n\rightarrow \infty $, for $t\in [0,T]$, $K_{T}^{p,+}-K_{T}^{+}%
\geq K_{t}^{p,+}-K_{t}^{+}\geq 0$. So as $n\rightarrow 0$,
\[
E\sup_{0\leq t\leq T}\left| K_{t}^{p,+}-K_{t}^{+}\right| ^{2}\leq E\left|
K_{T}^{p,+}-K_{T}^{+}\right| ^{2}\rightarrow 0,
\]
Similarly, we have $E\sup_{0\leq t\leq T}\left| K_{t}^{-}-K_{t}^{p,-}\right|
^{2}\leq E\left| K_{T}^{-}-K_{T}^{p,-}\right| ^{2}\rightarrow 0$.

It remains to check if $(Y_{t},Z_{t},K_{t})_{0\leq t\leq T}$ satisfies (3)
and (4) of the definition \ref{Def-tb}. Since $L_{t}\leq Y_{t}^{n}\leq U_{t}$%
, $0\leq t\leq T$, then letting $n\rightarrow \infty $, $L_{t}\leq Y_{t}\leq
U_{t}$, $0\leq t\leq T$, a.s.. Furthermore $(Y^{n},K^{n,+})$ tends to $%
(Y,K^{+})$ uniformly in $t$ in probability, as $n\rightarrow \infty $, then
the measure $dK^{n,+}\rightarrow dK^{+}$ weakly in probability, as $%
n\rightarrow \infty $, i.e. $\int_{0}^{T}(Y_{t}^{n}-L_{t})dK_{t}^{n,+}%
\rightarrow \int_{0}^{T}(Y_{t}-L_{t})dK_{t}^{+}$, in probability. While $%
L_{t}\leq Y_{t}\leq U_{t}$, $0\leq t\leq T$, so $%
\int_{0}^{T}(Y_{t}-L_{t})dK_{t}^{+}\geq 0$. On the other hand $%
\int_{0}^{T}(Y_{t}^{n}-L_{t})dK_{t}^{n,+}=0$, so $%
\int_{0}^{T}(Y_{t}-L_{t})dK_{t}^{+}=0$. Similarly, $%
\int_{0}^{T}(Y_{t}-U_{t})dK_{t}^{-}=0$, i.e. the triple $%
(Y_{t},Z_{t},K_{t})_{0\leq t\leq T}$ is the solution of RBSDE$(\xi ,f,L)$,
under the assumption (\ref{assup2-blow}). $\square $

\paragraph{\textbf{Step 5.}}

Now we consider a terminal condition $\xi \in \mathbf{L}^{2}(\mathcal{F}%
_{T}) $ and a coefficient $f$ which satisfies assumption 2.4\textbf{. }For $%
n\in \mathbf{N}$, set
\[
\xi _{n}=\xi \vee (-n),f_{n}(t,y)=f(t,y)-f(t,0)+f(t,0)\vee (-n).
\]
Obviously, $(\xi ^{n},f^{n})$ satisfies the assumptions of the step 4, and
since $\xi \in \mathbf{L}^{2}(\mathcal{F}_{T})$, $f(t,0)\in \mathbf{H}%
^{2}(0,T)$, then
\[
E[\left| \xi ^{n}-\xi \right| ^{2}]\rightarrow 0,E\int_{0}^{T}\left|
f(t,0)-f_{n}(t,0)\right| ^{2}\rightarrow 0,
\]
as $n\rightarrow \infty $.

From the results in step 3, for each $n\in \mathbf{N}$, there exists $%
(Y_{t}^{n},Z_{t}^{n},K_{t}^{n})_{0\leq t\leq T}\in \mathbf{S}^{2}(0,T)\times
\mathbf{H}_{d}^{2}(0,T)\times \mathbf{VF}^{2}(0,T)$, with $%
K^{n}=K^{n,+}-K^{n,-}$, which is the unique solution of the RBSDE$(\xi
^{n},f_{n},L,U)$. Like in step 4, we consider the solution $(\overline{Y},%
\overline{Z},\overline{K})$ of the one lower barrier RBSDE$(\xi ^{+},%
\overline{f},L)$, where $\xi ^{+}$ is the positive part of $\xi $, $%
\overline{f}(t,y)=f(t,y)-f(t,0)+(f(t,0))^{+}$, and the solution $(\underline{%
Y},\underline{Z},\underline{K})$ of the one super barrier RBSDE$(\xi ,f,U)$.
Then we can take the RBSDE$(\xi ^{+},\overline{f},L)$ (resp. RBSDE$(\xi
,f,U) $) as a RBSDE with two barriers associated to the parameters $(\xi
^{+},\overline{f},L,\overline{U})$ (resp. $(\xi ,f,\underline{L},U)$), where
$\overline{U}=\infty $ and $\underline{L}=-\infty $. Thanks to the
comparison theorem \ref{comprg}, we have that
\[
E[\sup_{0\leq t\leq T}\left| Y_{t}^{n}\right| ^{2}]\leq \max \{E[\sup_{0\leq
t\leq T}\left| \overline{Y}_{t}\right| ^{2},E[\sup_{0\leq t\leq T}\left|
\underline{Y}_{t}\right| ^{2}\}\leq C.
\]

For $n$, $p\in \mathbf{N}$, with $n\geq p$, we have $\xi ^{n}\leq
\xi ^{p}$ and $f_{n}(t,y)\leq f_{p}(t,y)$, $\forall (t,y)\in
[0,T]\times \mathbb{R}$. From approximations for $\xi ^{n}$, $\xi
^{p}$, $f_{n}(t,y)$ and $f_{p}(t,y)$ as following:
\begin{eqnarray*}
\xi ^{n,m} &:&=\xi ^{n}\wedge m,\xi ^{p,m}:=\xi ^{p}\wedge m \\
f_{n,m}(t,y) &=&f_{n}(t,y)-f_{n}(t,0)+f_{n}(t,0)\wedge
m=f(t,y)-f(t,0)+(f(t,0)\vee (-n))\wedge m, \\
f_{p,m}(t,y) &=&f_{p}(t,y)-f_{p}(t,0)+f_{p}(t,0)\wedge
m=f(t,y)-f(t,0)+(f(t,0)\vee (-p))\wedge m,
\end{eqnarray*}
then the parameters satisfy the assumptions in theorem \ref{compr4}, and
\[
\xi ^{n,m}\leq \xi ^{p,m},f_{n,m}(t,y)\leq f_{p,m}(t,y).
\]
Consider the solution $(Y^{n,m},Z^{n,m},K^{n,m})$(resp. $%
(Y^{p,m},Z^{p,m},K^{p,m})$) of the RBSDE$(\xi ^{n,m},f_{n,m},L,U)$(resp. $%
(\xi ^{p,m},f_{p,m},L,U)$); by the comparison theorem \ref{compr4}, for $%
0\leq s\leq t\leq T$, we have $K_{t}^{n,m,-}-K_{s}^{n,m,-}\leq
K_{t}^{p,m,-}-K_{s}^{p,m,-}$. Then by the convergence results in step 4, let
$m\rightarrow \infty $, we get
\[
K_{t}^{n,-}-K_{s}^{n,-}\leq K_{t}^{p,-}-K_{s}^{p,-},\mbox{ for
}n\geq p.
\]
So we have $0\leq K_{t}^{n,-}\leq K_{t}^{1,-}$, then $E[(K_{t}^{n,-})^{2}]%
\leq E[(K_{t}^{1,-})^{2}]\leq C$. By the same method as previous step, we
deduce that
\[
E[\int_{0}^{t}f(s,Y_{s}^{n})ds)^{2}]+E[(K_{T}^{n,+})^{2}]+E\int_{0}^{T}%
\left| Z_{s}^{n}\right| ^{2}ds\leq C.
\]
Now we are in the same situation as in step 4, and following the same
method, we get that the sequence $%
(Y_{t}^{n},Z_{t}^{n},K_{t}^{n,+},K_{t}^{n,-})$ converge to $%
(Y_{t},Z_{t},K_{t}^{+},K_{t}^{-})$ as $n\rightarrow \infty $, in $\mathbf{S}%
^{2}(0,T)\times \mathbf{H}_{d}^{2}(0,T)\times \mathbf{A}^{2}(0,T)\times
\mathbf{A}^{2}(0,T)$, and $(Y_{t},Z_{t},K_{t}^{+},K_{t}^{-})$ is the
solution to the RBSDE$(\xi ,f,L,U)$. $\square $

\section{Appendix}

\subsection{Proofs of Lemmas}

In this subsection, we prove lemma \ref{com-l} and lemma \ref{com-bu}, which
play important roles in previous section.

\textbf{Proof of Lemma \ref{com-l}:} Obviously $f(s,(L_{s})^{-})\in \mathbf{H%
}^{2}(0,T)$, in view of assumption 2.3'\textbf{. }Consider for $m$, $n\in
\mathbf{N}$, the following RBSDEs with one lower barrier,
\begin{eqnarray*}
\widetilde{Y}_{t}^{m} &=&\xi +\int_{t}^{T}f(s,(L_{s})^{-})ds-m\int_{t}^{T}(%
\widetilde{Y}_{s}^{m}-U_{s})^{+}ds+\widetilde{K}_{T}^{m,+}-\widetilde{K}%
_{t}^{m,+}-\int_{t}^{T}\widetilde{Z}_{s}^{m}dB_{s}, \\
\widetilde{Y}_{t}^{m} &\geq &L_{t},\int_{0}^{T}(\widetilde{Y}_{t}^{m}-L_{t})d%
\widetilde{K}_{t}^{m,+}=0,
\end{eqnarray*}
and
\begin{eqnarray*}
Y_{t}^{m,n} &=&\xi
+\int_{t}^{T}f(s,Y_{s}^{m,n})ds-m%
\int_{t}^{T}(Y_{s}^{m,n}-U_{s})^{+}ds+K_{T}^{m,n,+}-K_{t}^{m,n,+}-%
\int_{t}^{T}Z_{s}^{m,n}dB_{s}, \\
Y_{t}^{m,n} &\geq
&L_{t}^{n},\int_{0}^{T}(Y_{t}^{m,n}-L_{t}^{n})dK_{t}^{m,n,+}=0.
\end{eqnarray*}
Since $Y_{t}^{m,n}\geq L_{t}^{n}\geq (L_{t})^{-}$, we get $%
f(t,Y_{t}^{m,n})\leq f(t,L_{t}^{n})\leq f(t,(L_{t})^{-})$. Then for $m,n\in
\mathbf{N}$, $\forall t\in [0,T]$%
\[
f(t,Y_{t}^{m,n})-m(Y_{t}^{m,n}-U_{t})^{+}\leq
f(t,(L_{t})^{-})-m(Y_{t}^{m,n}-U_{t})^{+},L_{t}^{n}\leq L_{t}.
\]
By the general comparison theorem for RBSDE with one barrier theorem \ref
{com-one}, it follows $Y_{t}^{m,n}\leq \widetilde{Y}_{t}^{m}$, $\forall t\in
[0,T]$. Denote $\widetilde{K}_{t}^{m,-}=m\int_{0}^{t}(\widetilde{Y}%
_{s}^{m}-U_{s})^{+}ds$ and $K_{t}^{m,n,-}=m%
\int_{0}^{t}(Y_{s}^{m,n}-U_{s})^{+}ds$, then we get for $0\leq s\leq t\leq T$
\begin{equation}
K_{t}^{m,n,-}-K_{s}^{m,n,-}\leq \widetilde{K}_{t}^{m,-}-\widetilde{K}%
_{s}^{m,-}.  \label{com-lin1}
\end{equation}
Thanks to the convergence result in step 1 and in \cite{LS}, notice that $%
(L^{n})^{+}$ is bounded, we know that as $m\rightarrow \infty $, $\widetilde{%
K}_{t}^{m,-}\rightarrow \widetilde{K}_{t}^{-},K_{t}^{m,n,-}\rightarrow
K_{t}^{n,-}$, in $\mathbf{L}^{2}(\mathcal{F}_{t})$. Here $\widetilde{K}%
_{t}^{-}$ and $K_{t}^{n,-}$ are increasing processes with respect to the
upper barrier $U$ of the solution of the RBSDE$(\xi ,f(t,(L_{t})^{-}),L,U)$
and RBSDE$(\xi ,f,L^{n},U)$, respectively. Then from (\ref{com-lin1}), we
deduce that for $0\leq s\leq t\leq T,$
\[
K_{t}^{n,-}-K_{s}^{n,-}\leq \widetilde{K}_{t}^{-}-\widetilde{K}_{s}^{-}.
\]
It follows immediately that $K_{T}^{n,-}\leq \widetilde{K}_{T}^{-}$. $%
\square $

\textbf{Proof of Lemma \ref{com-bu}: }Obviously $f(s,U_{s}{}^{+})\in \mathbf{%
H}^{2}(0,T)$, in view of assumption 2.3'\textbf{-}(i)\textbf{.} Consider the
following RBSDEs with one barrier, for $n$, $m$, $p\in \mathbf{N}$, with $%
U^{n}=U\vee (-n)$,
\begin{eqnarray*}
\widetilde{Y}_{t}^{n,m,p} &=&\xi
+\int_{t}^{T}f(s,(U_{s})^{+})ds+m\int_{t}^{T}(\widetilde{Y}%
_{s}^{n,m,p}-L_{s}^{p})^{-}ds-(\widetilde{K}_{T}^{n,m,p,-}-\widetilde{K}%
_{t}^{n,m,p,-})-\int_{t}^{T}\widetilde{Z}_{s}^{n,m,p}dB_{s}, \\
\widetilde{Y}_{t}^{n,m,p} &\leq &U_{t}^{n},\int_{0}^{T}(\widetilde{Y}%
_{t}^{n,m,p}-U_{t})d\widetilde{K}_{t}^{n,m,p,-}=0,
\end{eqnarray*}
and
\begin{eqnarray*}
Y_{t}^{n,m,p} &=&\xi
+\int_{t}^{T}f(s,Y_{s}^{n,m,p})ds+m%
\int_{t}^{T}(Y_{s}^{n,m,p}-L_{s}^{p})^{-}ds-(K_{T}^{n,m,p,-}-K_{t}^{n,m,p,-})-\int_{t}^{T}Z_{s}^{n,m,p}dB_{s},
\\
Y_{t}^{n,m,p} &\leq
&U_{t}^{n},\int_{0}^{T}(Y_{t}^{n,m,p}-U_{t}^{n})dK_{t}^{n,m,p,-}=0.
\end{eqnarray*}
Since $Y_{t}^{n,m,p}\leq U_{t}^{n}\leq (U_{t})^{+}$, by monotonic property
of $f$, we get $f(t,Y_{t}^{n,m,p})\geq f(t,(U_{t})^{+})$. So
\[
f(t,(U_{t})^{+})+m(Y_{t}^{n,m,p}-L_{t}^{p})^{-}\leq
f(t,Y_{t}^{n,m,p})+m(Y_{t}^{n,m,p}-L_{t}^{p})^{-},U_{t}\leq U_{t}^{n},
\]
from general comparison theorem for RBSDE with one barrier \ref{com-one}, we
have $Y_{t}^{n,m,p}\geq \widetilde{Y}_{t}^{n,m,p}$. Set $K_{t}^{n,m,p,+}=m%
\int_{0}^{t}(Y_{s}^{n,m,p}-L_{s}^{p})^{-}ds$, $\widetilde{K}%
_{t}^{n,m,p,+}=m\int_{0}^{t}(\widetilde{Y}_{s}^{n,m,p}-L_{s}^{p})^{-}ds$,
then for $0\leq s\leq t\leq T$%
\[
K_{t}^{n,m,p,+}-K_{s}^{n,m,p,+}\leq \widetilde{K}_{t}^{n,m,p,+}-\widetilde{K}%
_{s}^{n,m,p,+}.
\]
Notice that $(L^{p})^{+}$ and $(U^{n})^{-}$ are bounded, by convergence
results in step 1, and the convergence result in \cite{LS}, as $m\rightarrow
\infty $, for $0\leq s\leq t\leq T$, we have
\[
K_{t}^{n,p,+}-K_{s}^{n,p,+}\leq \widetilde{K}_{t}^{n,p,+}-\widetilde{K}%
_{s}^{n,p,+},
\]
where $K_{t}^{n,p,+}$ and $\widetilde{K}_{t}^{n,p,+}$ are the increasing
processes corresponding to lower barrier $L^{p}$ for RBSDE$(\xi
,f,L^{p},U^{n})$ and RBSDE$(\xi ,f(t,(U_{t})^{+}),L^{p},U^{n})$,
respectively. Then thanks to the convergence result in step 2 for the
approximation of lower barrier $L$, we have that as $p\rightarrow \infty $,
\[
K_{t}^{n,p,+}\rightarrow K_{t}^{n,+}\mbox{ and }\widetilde{K}%
_{t}^{n,p,+}\rightarrow \widetilde{K}_{t}^{n,+}\mbox{ in }\mathbf{L}^{2}(%
\mathcal{F}_{t}).
\]
Here $K^{n,+}$ (resp. $\widetilde{K}^{n,+}$) is the increasing process
corresponding to lower barrier $L$ for RBSDE$(\xi ,f,L,U^{n})$ (resp. RBSDE$%
(\xi ,f(t,(U_{t})^{+}),L,U^{n})$). It follows for $0\leq s\leq t\leq T$%
\[
K_{t}^{n,+}-K_{s}^{n,+}\leq \widetilde{K}_{t}^{n,+}-\widetilde{K}_{s}^{n,+}.
\]
Finally by comparison theorem \ref{com-2k}, since $U_{t}\leq U_{t}^{n}$, $%
\forall t\in [0,T]$, we get
\[
\widetilde{K}_{t}^{+}-\widetilde{K}_{s}^{+}\geq \widetilde{K}_{t}^{n,+}-%
\widetilde{K}_{s}^{n,+},
\]
where $\widetilde{K}_{t}^{+}$ is the increasing process corresponding to
lower barrier $L$ for RBSDE$(\xi ,f(t,(U_{t})^{+}),L,U^{n})$. So for $0\leq
s\leq t\leq T$%
\[
K_{t}^{n,+}-K_{s}^{n,+}\leq \widetilde{K}_{t}^{+}-\widetilde{K}_{s}^{+}.
\]
Specially, $K_{T}^{n,+}\leq \widetilde{K}_{T}^{+}$. $\square $

\subsection{Comparison theorems}

First we need a general comparison theorem for the RBSDE with one lower
barrier.

\begin{theorem}[General case for RBSDE's]
\label{com-one}Suppose that the parameters $(\xi ^{1},f^{1},L^{1})$ and $%
(\xi ^{2},f^{2},L^{2})$ satisfy assumption 2.1-2.3-(i). Let the triples $%
(Y^{1},Z^{1},K^{1})$, $(Y^{2},Z^{2},K^{2})$ be respectively the solutions of
the RBSDE$(\xi ^{1},f^{1},L^{1})$ and RBSDE$(\xi ^{2},f^{2},L^{2})$, i.e.
\[
Y_{t}^{i}=\xi
^{i}+\int_{t}^{T}f^{i}(s,Y_{s}^{i},Z_{s}^{i})ds+K_{T}^{i}-K_{t}^{i}-%
\int_{t}^{T}Z_{s}^{i}dB_{s},
\]
$Y_{t}^{i}\geq L_{t}^{i},\;0\leq t\leq T$, and $%
\int_{0}^{T}(Y_{s}^{i}-L_{s}^{i})dK_{s}^{i}=0$, $i=1,2$. Assume in addition
the following: $\forall t\in [0,T],$%
\begin{equation}
\mbox{ }\xi ^{1}\leq \xi ^{2},f^{1}(t,Y_{t}^{1},Z_{t}^{1})\leq
f^{2}(t,Y_{t}^{1},Z_{t}^{1})\mathbf{,}L_{t}^{1}\leq L_{t}^{2},
\label{con-compg1}
\end{equation}
then $Y_{t}^{1}\leq Y_{t}^{2}$, for $t\in [0,T]$.
\end{theorem}

\textbf{Proof. }Applying It\^{o}'s formula to $[(Y^{1}-Y^{2})^{+}]^{2}$ on
interval $[t,T]$, and taking expectation on the both sides, since on the set
$\{Y_{t}^{1}>Y_{t}^{2}\}$, $Y_{t}^{1}>Y_{t}^{2}\geq L_{t}^{2}\geq L_{t}^{1}$%
, we have
\[
\int_{t}^{T}(Y_{s}^{1}-Y_{s}^{2})^{+}d(K_{s}^{1}-K_{s}^{2})=-%
\int_{t}^{T}(Y_{s}^{1}-Y_{s}^{2})^{+}dK_{s}^{2}\leq 0,
\]
then we get immediately
\begin{eqnarray*}
&&E[(Y_{t}^{1}-Y_{t}^{2})^{+}]^{2}+E\int_{t}^{T}1_{\{Y_{t}^{1}>Y_{t}^{2}\}}%
\left| Z_{s}^{1}-Z_{s}^{2}\right| ^{2}ds \\
&\leq &\frac{1}{2}E\int_{t}^{T}1_{\{Y_{t}^{1}>Y_{t}^{2}\}}\left|
Z_{s}^{1}-Z_{s}^{2}\right| ^{2}ds+(2\mu
+4k^{2})E\int_{t}^{T}[(Y_{s}^{1}-Y_{s}^{2})^{+}]^{2}ds,
\end{eqnarray*}
in view of (\ref{con-compg1}) and the Lipschitz condition and monotonic
condition of $f^{2}$. Hence
\[
E[(Y_{t}^{1}-Y_{t}^{2})^{+}]^{2}\leq (2\mu
+4K^{2})E\int_{t}^{T}[(Y_{s}^{1}-Y_{s}^{2})^{+}]^{2}ds,
\]
from Gronwall's inequality, we deduce $(Y_{t}^{1}-Y_{t}^{2})^{+}=0$, $0\leq
t\leq T$.$\square $

Then we prove a comparison theorem for the increasing processes under
Lipschitz assumption on $f$ via the penalization method in \cite{LS}.

\begin{theorem}
\label{com-2k}Suppose that the parameters $(\xi ^{1},f^{1},L^{1},U^{1})$ and
$(\xi ^{2},f^{2},L^{2},U^{2})$ satisfy the following conditions: for $i=1,2$,

(i) $\xi ^{i}\in \mathbf{L}^{2}(\mathcal{F}_{T})$;

(ii) $f^{i}$ satisfy assumption 2.2-(i), (iii), (vi) and a Lipschitz
condition in $(y,z)$ uniformly in $(t,\omega )$, i.e. there exists a
constant $k$ such that, for $y,y^{\prime }\in \mathbb{R}$,
$z,z^{\prime }\in \mathbb{R}^{d}$,
\[
\left| f^{i}(t,y,z)-f^{i}(t,y^{\prime },z^{\prime })\right| \leq k(\left|
y-y^{\prime }\right| +\left| z-z^{\prime }\right| );
\]

(iii) $L^{i}$ and $U^{i}$ are real-valued, $\mathcal{F}_{t}$-adapted,
continuous with $(L^{i})^{+}$, $(U^{i})^{-}\in \mathbf{S}^{2}(0,T)$.

Let $(Y^{i},Z^{i},K^{i,+},K^{i,-})$ be the solution of the RBSDE$(\xi
^{i},f^{i},L^{i},U^{i})$, i.e.
\[
Y_{t}^{i}=\xi
^{i}+%
\int_{t}^{T}f^{i}(s,Y_{s}^{i},Z_{s}^{i})ds+K_{T}^{i,+}-K_{t}^{i,+}-(K_{T}^{i,-}-K_{t}^{i,-})-\int_{t}^{T}Z_{s}^{i}dB_{s},
\]
$Y_{t}^{i}\geq L_{t}^{i},\;0\leq t\leq T$, and $%
\int_{0}^{T}(Y_{s}^{i}-L_{s}^{i})dK_{s}^{i,+}=%
\int_{0}^{T}(Y_{s}^{i}-U_{s}^{i})dK_{s}^{i,-}=0$,

Moreover, we assume $\forall (t,y,z)\in [0,T]\times \mathbb{R\times
R}^{d}$,
\begin{equation}
\mbox{ }\xi ^{1}\leq \xi ^{2},f^{1}(t,y,z)\leq f^{2}(t,y,z),
\nonumber
\end{equation}
Then we have: for $0\leq t\leq T$,

(i) If $L^{1}=L^{2}$, $U^{1}=U^{2}$, then $Y_{t}^{1}\leq Y_{t}^{2}$, $%
K_{t}^{1,+}\geq K_{t}^{2,+}$, $K_{t}^{1,-}\leq K_{t}^{2,-}$, and $%
dK^{1,+}\geq dK^{2,+}$, $dK^{1,-}\leq dK^{2,-}$;

(ii) If $L_{t}^{1}\leq L_{t}^{2}$, $U_{t}^{1}=U_{t}^{2}$, then $%
Y_{t}^{1}\leq Y_{t}^{2}$, $K_{t}^{1,-}\leq K_{t}^{2,-}$, and $dK^{1,-}\leq
dK^{2,-}$;

(iii) If $L_{t}^{1}=L_{t}^{2}$, $U_{t}^{1}\leq U_{t}^{2}$, then $%
Y_{t}^{1}\leq Y_{t}^{2}$, $K_{t}^{1,+}\geq K_{t}^{2,+}$, and $dK^{1,+}\geq
dK^{2,+}$.
\end{theorem}

\textbf{Proof. }(i) Set $L:=L^{1}=L^{2}$, $U:=U^{1}=U^{2}$, and consider the
penalization equations for $m$, $n\in \mathbf{N}$, $i=1,2$%
\[
Y_{t}^{m,n,i}=\xi
^{i}+\int_{t}^{T}f^{i}(s,Y_{s}^{m,n,i},Z_{s}^{m,n,i})ds+m%
\int_{t}^{T}(Y_{s}^{m,n,i}-L_{s})^{-}ds-n%
\int_{t}^{T}(Y_{s}^{m,n,i}-U_{s})^{+}ds-\int_{t}^{T}Z_{s}^{m,n,i}dB_{s}.
\]
By comparison theorem for BSDEs, since $\xi ^{1}\leq \xi ^{2}$ and
\[
f^{1}(t,y,z)+m(y-L_{t})^{-}-n(y-U_{t})^{+}\leq
f^{2}(t,y,z)+m(y-L_{t})^{-}-n(y-U_{t})^{+},
\]
we have $Y_{t}^{m,n,1}\leq Y_{t}^{m,n,2}$, $\forall t\in [0,T]$. Denote $%
K_{t}^{m,n,i,+}=m\int_{0}^{t}(Y_{s}^{m,n,i}-L_{s})^{-}ds$, then for $0\leq
s\leq t\leq T$,
\[
K_{t}^{m,n,1,+}-K_{s}^{m,n,1,+}\geq K_{t}^{m,n,2,+}-K_{s}^{m,n,2,+}.
\]
From the convergence results in \cite{LS}, which also holds for Lipschitz
function,
\[
\lim_{n\rightarrow \infty }\lim_{m\rightarrow \infty
}Y_{t}^{m,n,i}=Y_{t}^{i},\lim_{n\rightarrow \infty }\lim_{m\rightarrow
\infty }K_{t}^{m,n,i,+}=K_{t}^{i,+},
\]
in $\mathbf{L}^{2}(\mathcal{F}_{t})$, where $Y^{i}$, $K^{i,+}$, $K^{i,-}$
are elements of the solution of RBSDE$(\xi ^{i},f^{i},L,U)$. Consequently,
for $0\leq s\leq t\leq T$,
\[
Y_{t}^{1}\leq Y_{t}^{2},K_{t}^{1,+}-K_{s}^{1,+}\geq K_{t}^{2,+}-K_{s}^{2,+};
\]
if we especially set $s=0$, we get $K_{t}^{1,+}\geq K_{t}^{2,+}$. Similarly $%
K_{t}^{1,-}\leq K_{t}^{2,-}$.

(ii) Set $U:=U^{1}=U^{2}$, we consider the penalized reflected BSDE's, for $%
n\in \mathbf{N}$, $i=1,2$,
\begin{eqnarray*}
Y_{t}^{n,i} &=&\xi
^{i}+%
\int_{t}^{T}f^{i}(s,Y_{s}^{n,i},Z_{s}^{n,i})ds+K_{T}^{n,+}-K_{t}^{n,+}-n%
\int_{t}^{T}(Y_{s}^{n,i}-U_{s})^{+}ds-\int_{t}^{T}Z_{s}^{n,i}dB_{s}, \\
Y_{t}^{n,i} &\geq
&L_{t}^{i},\int_{0}^{T}(Y_{t}^{n,i}-L_{t}^{i})dK_{t}^{n,+}=0.
\end{eqnarray*}
Since $\forall t\in [0,T]$,
\[
\xi ^{1}\leq \xi ^{2},\,\,\,\,f^{1}(t,y,z)-n(y-U_{t})^{+}\leq
f^{2}(t,y,z)-n(y-U_{t})^{+},\,\,\,\,\,L_{t}^{1}\leq L_{t}^{2}
\]
by the comparison theorem for RBSDE with one barrier, we have $%
Y_{t}^{n,1}\leq Y_{t}^{n,2}$. Let $K_{t}^{n,i,-}=n%
\int_{0}^{t}(Y_{s}^{n,i}-U_{s})^{+}ds$, then for $0\leq s\leq t\leq T$,
\[
K_{t}^{n,1,-}-K_{s}^{n,1,-}\leq K_{t}^{n,2,-}-K_{s}^{n,2,-}.
\]
Thanks to the convergence result in \cite{LS}, which still works for
Lipschitz functions, we have for $0\leq s\leq t\leq T$,
\[
Y_{t}^{1}\leq Y_{t}^{2},K_{t}^{1,-}-K_{s}^{1,-}\leq K_{t}^{2,-}-K_{s}^{2,-};
\]
if we especially set $s=0$, we get $K_{t}^{1,+}\geq K_{t}^{2,+}$, $%
K_{t}^{1,-}\leq K_{t}^{2,-}$.

(iii) The proof is in the same as (ii), so we omit it. $\square $

We next prove a comparison theorem for RBSDE with two barriers in a general
case.

\begin{theorem}[General case for RBSDE's]
\label{comprg}Suppose that the parameters $(\xi ^{1},f^{1},L^{1},U^{1})$ and
$(\xi ^{2},f^{2},L^{2},U^{2})$ satisfy assumption 2.1\textbf{, }2.2 and 2.3%
\textbf{.} Let $(Y^{1},Z^{1},K^{1,+},K^{1,-})$, $%
(Y^{2},Z^{2},K^{2,+},K^{2,-})$ be respectively the solutions of the RBSDE$%
(\xi ^{1},f^{1},L^{1},U^{1})$ and RBSDE$(\xi ^{2},f^{2},L^{2},U^{2})$ as
definition \ref{Def-tb}. Assume in addition the following: $\forall t\in
[0,T]$%
\begin{eqnarray}
\mbox{ }\xi ^{1} &\leq &\xi
^{2},\,\,\,\,\,f^{1}(t,Y_{t}^{1},Z_{t}^{1})\leq
f^{2}(t,Y_{t}^{1},Z_{t}^{1}),  \label{con-compg} \\
L_{t}^{1} &\leq &L_{t}^{2},\,\,\,\,\,U_{t}^{1}\leq U_{t}^{2},  \nonumber
\end{eqnarray}
then $Y_{t}^{1}\leq Y_{t}^{2}$, for $t\in [0,T]$.
\end{theorem}

\textbf{Proof. }Applying Ito's formula to $[(Y^{1}-Y^{2})^{+}]^{2}$ on the
interval $[t,T]$, and taking expectation on the both sides, we get
immediately
\begin{eqnarray*}
&&E[(Y_{t}^{1}-Y_{t}^{2})^{+}]^{2}+E\int_{t}^{T}1_{\{Y_{t}^{1}>Y_{t}^{2}\}}%
\left| Z_{s}^{1}-Z_{s}^{2}\right| ^{2}ds \\
&\leq
&2E\int_{t}^{T}1_{\{Y_{s}^{1}>Y_{s}^{2}%
\}}(Y_{s}^{1}-Y_{s}^{2})(f^{2}(s,Y_{s}^{1},Z_{s}^{1})-f^{2}(s,Y_{s}^{2},Z_{s}^{2}))ds
\\
&\leq &\frac{1}{2}E\int_{t}^{T}1_{\{Y_{t}^{1}>Y_{t}^{2}\}}\left|
Z_{s}^{1}-Z_{s}^{2}\right| ^{2}ds+(2\mu
+4k^{2})E\int_{t}^{T}[(Y_{s}^{1}-Y_{s}^{2})^{+}]^{2}ds,
\end{eqnarray*}
in view of (\ref{con-compg}) and the Lipschitz condition and monotonic
condition on $f^{2}$, and the fact that
\[
\int_{t}^{T}(Y_{s}^{1}-Y_{s}^{2})^{+}d(K_{s}^{1,+}-K_{s}^{2,+})\leq
0,\int_{t}^{T}(Y_{s}^{1}-Y_{s}^{2})^{+}d(K_{s}^{1,-}-K_{s}^{2,-})\geq 0,
\]
which is similar to reflected BSDE with one barrier. Hence
\[
E[(Y_{t}^{1}-Y_{t}^{2})^{+}]^{2}\leq (2\mu
+4k^{2})E\int_{t}^{T}[(Y_{s}^{1}-Y_{s}^{2})^{+}]^{2}ds,
\]
and from Gronwall's inequality, we deduce $(Y_{t}^{1}-Y_{t}^{2})^{+}=0$, $%
0\leq t\leq T$.$\square $

From the convergence of penalization equations, we get the following
comparison theorem.

\begin{theorem}[Special case]
\label{compr1}Suppose that $f^{1}(s,y)$, $f^{2}(s,y)$ satisfy assumption
2.4, and$\ \xi ^{i}$, $f^{i}(\cdot ,0)$, $L$, $U$, $i=1,2$ satisfies the
bounded assumption (\ref{assup-bon}). The two triples $(Y^{1},Z^{1},K^{1})$,
$(Y^{2},Z^{2},K^{2})$ are respectively the solutions of the RBSDE$(\xi
^{1},f^{1},L,U)$ and RBSDE$(\xi ^{2},f^{2},L,U)$ as definition \ref{Def-tb}.
If we have
\[
\mbox{ }\xi ^{1}\leq \xi ^{2},\mbox{ and }f^{1}(t,y)\leq
f^{2}(t,y),\forall (t,y)\in [0,T]\times \mathbb{R,}
\]
then $Y_{t}^{1}\leq Y_{t}^{2}$, $K_{t}^{1,+}\geq K_{t}^{2,+}$, $%
K_{t}^{1,-}\leq K_{t}^{2,-}$, for $t\in [0,T]$, and $dK^{1,+}\geq dK^{2,+}$,
$dK^{1,-}\leq dK^{2,-}$.
\end{theorem}

\textbf{Proof. }We consider the penalized equations relative to the RBSDE$%
(\xi ^{i},f^{i},L,U)$, for $i=1,2$, $n\in \mathbf{N},$%
\[
Y_{t}^{m,n,i}=\xi
^{i}+\int_{t}^{T}f^{i}(s,Y_{s}^{m,n,i})ds+n%
\int_{t}^{T}(Y_{s}^{m,n,i}-L_{s})^{-}ds-m%
\int_{t}^{T}(Y_{s}^{m,n,i}-U_{s})^{+}ds-\int_{t}^{T}Z_{s}^{m,n,i}dB_{s}.
\]
For each $m$, $n\in \mathbf{N}$,
\[
f^{m,n,1}(s,y)=f^{1}(s,y)+n(y-L_{s})^{-}-m(y-U_{s})^{+}\leq
f^{m,n,2}(s,y)=f^{2}(s,y)+n(y-L_{s})^{-}-m(y-U_{s})^{+},
\]
and $\xi ^{1}\leq \xi ^{2}$. So by the comparison theorem in \cite{Pardoux99}%
, we get
\[
Y_{t}^{m,n,1}\leq Y_{t}^{m,n,2},0\leq t\leq T.
\]

Since $K_{t}^{m,n,i,+}=n\int_{0}^{t}(Y_{s}^{m,n,i}-L_{s})^{-}ds$, then we
deduce, for $0\leq s\leq t\leq T$,
\[
K_{t}^{m,n,1,+}-K_{s}^{m,n,1,+}\geq K_{t}^{m,n,2,+}-K_{s}^{m,n,2,+},
\]

By the convergence results of the step1, we know tha the inequalities hold
for $0\leq s\leq t\leq T$:
\[
Y_{t}^{1}\leq Y_{t}^{2},K_{t}^{1,+}-K_{s}^{1,+}\geq K_{t}^{2,+}-K_{s}^{2,+},
\]
Particularly, set $s=0$, we get $K_{t}^{1,+}\geq K_{t}^{2,+}$.
Symmetrically, $K_{t}^{1,-}\leq K_{t}^{2,-}.\square $

\begin{corollary}
\label{compr2}Suppose that $f^{1}(s,y)$, $f^{2}(s,y)$ satisfy assumption
2.4, and$\ \xi ^{i}$, $f^{i}(\cdot ,0)$, $L^{i}$, $U^{i}$, $i=1,2$ satisfies
(\ref{assup-bon}). The two triples $(Y^{1},Z^{1},K^{1})$, $%
(Y^{2},Z^{2},K^{2})$ are respectively the solutions of the RBSDE$(\xi
^{1},f^{1},L^{1},U^{1})$ and RBSDE$(\xi ^{2},f^{2},L^{2},U^{2})$, with $%
K^{i}=K^{i,+}-K^{i,-}$, $i=1,2$. In addition, we assume
\[
\xi ^{1}\leq \xi ^{2},\mbox{ and }f^{1}(t,y)\leq f^{2}(t,y),\forall
(t,y)\in [0,T]\times \mathbb{R,}
\]
then for $0\leq t\leq T$, (i) If $L_{t}^{1}\leq L_{t}^{2}$, $%
U_{t}^{1}=U_{t}^{2}$, then $Y_{t}^{1}\leq Y_{t}^{2}$, $K_{t}^{1,-}\leq
K_{t}^{2,-}$, and $dK^{1,-}\leq dK^{2,-}$;

(ii) If $L_{t}^{1}=L_{t}^{2}$, $U_{t}^{1}\leq U_{t}^{2}$, then $%
Y_{t}^{1}\leq Y_{t}^{2}$, $K_{t}^{1,+}\geq K_{t}^{2,+}$, and $dK^{1,+}\geq
dK^{2,+}$.
\end{corollary}

\textbf{Proof. }(i) To simplify symbols, we denote $U=U^{1}=U^{2}$. For $%
n\in \mathbf{N}$, we consider the following RBSDE with one barrier $L^{i}$, $%
i=1,2.$%
\begin{eqnarray*}
Y_{t}^{n,i} &=&\xi
^{i}+\int_{t}^{T}f^{i}(s,Y_{s}^{n,i})ds+K_{T}^{n,i,+}-K_{t}^{n,i,+}-n%
\int_{t}^{T}(Y_{s}^{n,i}-U_{s})^{+}ds-\int_{t}^{T}Z_{s}^{n,i}dB_{s}, \\
Y_{t}^{n,i} &\geq &L_{t}^{i},\int_{0}^{T}(Y_{t}^{n,i}-L_{t}^{i})dt=0.
\end{eqnarray*}
Since $\xi ^{1}\leq \xi ^{2}$, $f^{1}(t,y)\leq f^{2}(t,y)$, $L_{t}^{1}\leq
L_{t}^{2}$, by general comparison theorem of RBSDE with one barrier, we know
$Y_{t}^{n,1}\leq Y_{t}^{n,2}$. Denote $K_{t}^{n,1,-}=n%
\int_{0}^{t}(Y_{s}^{n,1}-U_{s})^{+}ds$, $K_{t}^{n,2,-}=n%
\int_{0}^{t}(Y_{s}^{n,2}-U_{s})^{+}ds$, then for $0\leq s\leq t\leq T$,
\[
K_{t}^{n,1,-}-K_{s}^{n,1,-}\leq K_{t}^{n,2,-}-K_{s}^{n,2,-}.
\]
Thanks to the convergence result of step 1of theorem \ref{exist22}, it
follows immediately that for $0\leq s\leq t\leq T$,
\[
Y_{t}^{1}\leq Y_{t}^{2}\mbox{ and }K_{t}^{1,-}-K_{s}^{1,-}\leq
K_{t}^{2,-}-K_{s}^{2,-}.
\]
Especially with $s=0$, we get $K_{t}^{1,-}\leq K_{t}^{2,-}$.

(ii) follows similarly as (i), so we omit it.$\square $

\begin{theorem}
\label{compr3}Suppose that $f^{1}(s,y)$, $f^{2}(s,y)$ satisfy assumption 2.4,%
$\ \xi ^{i}$, $f^{i}(\cdot ,0)$, $U^{i}$, for $i=1,2$ satisfies (\ref
{assup-bon2}), and $L^{i}$ satisfies assumption 2.3'. The two groups $%
(Y^{1},Z^{1},K^{1})$, $(Y^{2},Z^{2},K^{2})$ are respectively the solutions
of the RBSDE$(\xi ^{1},f^{1},L^{1},U^{1})$ and RBSDE$(\xi
^{2},f^{2},L^{2},U^{2})$, as definition \ref{Def-tb}. Moreover, assume
\[
\mbox{ }\xi ^{1}\leq \xi ^{2},f^{1}(t,y)\leq f^{2}(t,y),\forall
(t,y)\in [0,T]\times \mathbb{R,}
\]
then for $0\leq t\leq T$,

(i) If $L^{1}=L^{2}$ and $U^{1}=U^{2}$, then $Y_{t}^{1}\leq Y_{t}^{2}$, $%
K_{t}^{1,+}\geq K_{t}^{2,+}$, $K_{t}^{1,-}\leq K_{t}^{2,-}$, and $%
dK^{1,+}\geq dK^{2,+}$, $dK^{1,-}\leq dK^{2,-}$;

(ii) If $L_{t}^{1}\leq L_{t}^{2}$, $U_{t}^{1}=U_{t}^{2}$, then $%
Y_{t}^{1}\leq Y_{t}^{2}$, $K_{t}^{1,-}\leq K_{t}^{2,-}$, and $dK^{1,-}\leq
dK^{2,-}$;

(iii) If $L_{t}^{1}=L_{t}^{2}$, $U_{t}^{1}\leq U_{t}^{2}$, then $%
Y_{t}^{1}\leq Y_{t}^{2}$, $K_{t}^{1,+}\geq K_{t}^{2,+}$, and $dK^{1,+}\geq
dK^{2,+}$.
\end{theorem}

\textbf{Proof. }Like in the step 2 of the proof of theorem \ref{exist22}, we
approximate the barrier $L^{i}$ by super bounded barrier $L^{n,i}$, where $%
L^{n,i}=L^{i}\wedge n$.

(i) Set $L:=L^{1}=L^{2}$, $U:=U^{1}=U^{2}$, and $L^{n}=L\wedge n$. Then
consider the RBSDE$(\xi ^{i},f^{i},L^{n},U)$, for $i=1,2,$%
\begin{eqnarray*}
Y_{t}^{n,i} &=&\xi
^{i}+%
\int_{t}^{T}f^{i}(s,Y_{s}^{n,i})ds+K_{T}^{n,i,+}-K_{t}^{n,i,+}-(K_{T}^{n,i,-}-K_{t}^{n,i,-})-\int_{t}^{T}Z_{s}^{n,i}dB_{s},
\\
L_{t}^{n} &\leq &Y_{t}^{n,i}\leq
U_{t},\int_{0}^{T}(Y_{s}^{n,i}-L_{s})dK_{s}^{n,i,+}=%
\int_{0}^{T}(Y_{s}^{n,i}-U_{s})dK_{s}^{n,i,-}=0.
\end{eqnarray*}
Since
\[
\mbox{ }\xi ^{1}\leq \xi ^{2},f^{1}(t,y)\leq f^{2}(t,y),
\]
from comparison theorem \ref{compr1}, we have for $0\leq s\leq t\leq T$,
\[
Y_{t}^{n,1}\leq Y_{t}^{n,2},K_{t}^{n,1,+}-K_{s}^{n,1,+}\geq
K_{t}^{n,2,+}-K_{s}^{n,2,+}.
\]
Thanks to the convergence results in step 2 of the proof for theorem \ref
{exist22}, we get that
\[
Y_{t}^{1}\leq Y_{t}^{2},K_{t}^{1,+}-K_{s}^{1,+}\geq K_{t}^{2,+}-K_{s}^{2,+},%
\mbox{ for }0\leq s\leq t\leq T.
\]
Especially, with $s=0$, we get $K_{t}^{1,+}\geq K_{t}^{2,+}$. Similarly $%
K_{t}^{1,-}\leq K_{t}^{2,-}$, for $t\in [0,T]$.

(ii) Set $U:=U^{1}=U^{2}$, and $L^{n,i}=L^{i}\wedge n$. Then we consider the
solutions $(Y^{n,i},Z^{n,i},K^{n,i})$ of the RBSDEs $(\xi
^{i},f^{i},L^{n,i},U)$, for $i=1,2$. Since
\[
\mbox{ }\xi ^{1}\leq \xi ^{2},f^{1}(t,y)\leq
f^{2}(t,y),L_{t}^{n,1}\leq L_{t}^{n,2},
\]
from corollary \ref{compr2}, we have for $0\leq s\leq t\leq T$, $%
Y_{t}^{n,1}\leq Y_{t}^{n,2}$ and $K_{t}^{n,1,-}-K_{s}^{n,1,-}\leq
K_{t}^{n,2,-}-K_{s}^{n,2,-}$. Then by the convergence results in step 2 of
the proof for theorem \ref{exist22}, it follows that
\[
Y_{t}^{1}\leq Y_{t}^{2},K_{t}^{1,-}-K_{s}^{1,-}\leq K_{t}^{2,-}-K_{s}^{2,-}%
\mbox{, for }0\leq s\leq t\leq T,.
\]
Especially, with $s=0$, we get $K_{t}^{1,-}\leq K_{t}^{2,-}$, for $t\in
[0,T] $.

(iii) The proof is similar to (ii), which follows from corollary \ref{compr2}
and by the convergence results in step 2 of the proof for theorem \ref
{exist22}, so we omit it. $\square $

\begin{theorem}
\label{compr4}Suppose that $f^{1}(s,y)$, $f^{2}(s,y)$ satisfy assumption 2.4,%
$\ \xi ^{i}$, $f^{i}(\cdot ,0)$, for $i=1,2$ satisfies the bounded
assumption (\ref{assup-bon3}), and $L^{i}$ and $U^{i}$satisfy assumption
2.3'. The two groups $(Y^{1},Z^{1},K^{1})$, $(Y^{2},Z^{2},K^{2})$ are
respectively the solutions of the RBSDE$(\xi ^{1},f^{1},L^{1},U^{1})$ and
RBSDE$(\xi ^{2},f^{2},L^{2},U^{2})$. Moreover, assume
\[
\xi ^{1}\leq \xi ^{2},f^{1}(t,y)\leq f^{2}(t,y),\forall (t,y)\in
[0,T]\times \mathbb{R,}
\]
then for $t\in [0,T]$,

(i) If $L^{1}=L^{2}$ and $U^{1}=U^{2}$, then $Y_{t}^{1}\leq Y_{t}^{2}$, $%
K_{t}^{1,+}\geq K_{t}^{2,+}$, $K_{t}^{1,-}\leq K_{t}^{2,-}$, and $%
dK^{1,+}\geq dK^{2,+}$, $dK^{1,-}\leq dK^{2,-}$;

(ii) If $L_{t}^{1}\leq L_{t}^{2}$, $U_{t}^{1}=U_{t}^{2}$, then $%
Y_{t}^{1}\leq Y_{t}^{2}$, $K_{t}^{1,-}\leq K_{t}^{2,-}$, and $dK^{1,-}\leq
dK^{2,-}$;

(iii) If $L_{t}^{1}=L_{t}^{2}$, $U_{t}^{1}\leq U_{t}^{2}$, then $%
Y_{t}^{1}\leq Y_{t}^{2}$, $K_{t}^{1,+}\geq K_{t}^{2,+}$, and $dK^{1,+}\geq
dK^{2,+}$.
\end{theorem}

\textbf{Proof. }Like in theorem \ref{compr3}, we approximate the barrier $U$
by lower bounded barrier $U^{n}$, where $U^{n}=U\vee (-n)$, then the results
follow from the comparison theorem \ref{compr3} and the convergence results
of step 3 in the proof of theorem \ref{exist22}, so we omit it.$\square $

\begin{theorem}
Suppose that for $i=1,2$,$\ \xi ^{i}$satisfies assumption 2.1, $f^{i}$ does
not depends on $z$ and satisfies assumption 2.2, $L^{i}$ and $U^{i}$satisfy
assumption 2.3. The two triples $(Y^{1},Z^{1},K^{1,+},K^{1,-})$, $%
(Y^{2},Z^{2},K^{2,+},K^{2,-})$ are the solutions of the RBSDE$(\xi
^{1},f^{1},L^{1},U^{1})$ and RBSDE$(\xi ^{2},f^{2},L^{2},U^{2}),$
respectively. Moreover, assume for $(t,y)\in [0,T]\times \mathbb{R,}$%
\[
\mbox{ }\xi ^{1}\leq \xi ^{2},f^{1}(t,y)\leq f^{2}(t,y),\mbox{ and }%
f^{1}(t,0)=f^{2}(t,0),
\]
then for $t\in [0,T]$,

(i) If $L^{1}=L^{2}$ and $U^{1}=U^{2}$, then $Y_{t}^{1}\leq Y_{t}^{2}$, $%
K_{t}^{1,+}\geq K_{t}^{2,+}$, $K_{t}^{1,-}\leq K_{t}^{2,-}$, and $%
dK^{1,+}\geq dK^{2,+}$, $dK^{1,-}\leq dK^{2,-}$;

(ii) If $L_{t}^{1}\leq L_{t}^{2}$, $U_{t}^{1}=U_{t}^{2}$, then $%
Y_{t}^{1}\leq Y_{t}^{2}$, $K_{t}^{1,-}\leq K_{t}^{2,-}$, and for $%
dK^{1,-}\leq dK^{2,-}$;

(iii) If $L_{t}^{1}=L_{t}^{2}$, $U_{t}^{1}\leq U_{t}^{2}$, then $%
Y_{t}^{1}\leq Y_{t}^{2}$, $K_{t}^{1,+}\geq K_{t}^{2,+}$, and for $%
dK^{1,+}\geq dK^{2,+}$.
\end{theorem}

\textbf{Proof. }(i)Set $L:=L^{1}=L^{2}$, $U:=U^{1}=U^{2}$. Like in the proof
of the theorem \ref{exist22}, for $i=1,2$, set
\[
(\overline{Y}_{t}^{i},\overline{Z}_{t}^{i},\overline{K}_{t}^{i,+},\overline{K%
}_{t}^{i,-}):=(e^{\lambda t}Y_{t}^{i},e^{\lambda
t}Z_{t}^{i},\int_{0}^{t}e^{\lambda s}dK_{s}^{i,+},\int_{0}^{t}e^{\lambda
s}dK_{s}^{i,-}).
\]
Then it's easy to check that for $i=1,2$, $(\overline{Y}_{t}^{i},\overline{Z}%
_{t}^{i},\overline{K}_{t}^{i,+},\overline{K}_{t}^{i,-})_{0\leq t\leq T}$ is
the solution of the RBSDE$(\overline{\xi }^{i},\overline{f}^{i},\overline{L},%
\overline{U})$, where
\[
(\overline{\xi }^{i},\overline{f}^{i}(t,y),\overline{L}_{t},\overline{U}%
_{t})=(e^{\lambda T}\xi ^{i},e^{\lambda t}f^{i}(t,e^{-\lambda t}y)-\lambda
y,e^{\lambda t}L_{t},e^{\lambda t}U_{t}).
\]
If we assume $\lambda =\mu $, then $(\overline{\xi }^{i},\overline{f}^{i},%
\overline{L},\overline{U})$ satisfies assumption 2.1, 2.4 and 2.3'. Since
the transform keeps the monotonicity, the results are equivalent to
\begin{equation}
\overline{Y}_{t}^{1}\leq \overline{Y}_{t}^{2},\overline{K}_{t}^{1,+}-%
\overline{K}_{s}^{1,+}\geq \overline{K}_{t}^{2,+}-\overline{K}_{s}^{2,+},%
\overline{K}_{t}^{1,-}-\overline{K}_{s}^{1,-}\geq \overline{K}_{t}^{2,-}-%
\overline{K}_{s}^{2,-},  \label{comp-rem}
\end{equation}
for $0\leq s\leq t\leq T$. Then we make the approximations
\begin{eqnarray*}
\overline{\xi }^{m,n,i} &:&=\overline{\xi }^{n,i}\wedge m:=(\overline{\xi }%
^{i}\vee (-n))\wedge m \\
\overline{f}_{m,n}^{i}(t,y) &:&=\overline{f}_{n}^{i}(t,y)-\overline{f}%
_{n}^{i}(t,0)+\overline{f}_{n}^{i}(t,0)\wedge m \\
&:&=\overline{f}^{i}(t,y)-\overline{f}^{i}(t,0)+(\overline{f}^{i}(t,0)\vee
(-n))\wedge m.
\end{eqnarray*}
Let for $i=1,2$, $(\overline{Y}_{t}^{m,n,i},\overline{Z}_{t}^{m,n,i},%
\overline{K}_{t}^{m,n,i,+},\overline{K}_{t}^{m,n,i,-})_{0\leq t\leq T}$ be
the solution of the RBSDE $(\overline{\xi }^{m,n,i},\overline{f}_{m,n}^{i},%
\overline{L},\overline{U})$; then $\overline{\xi }^{m,n,i}$, $\overline{f}%
_{m,n}^{i}$ satisfy
\[
\left| \overline{\xi }^{m,n,i}\right| +\sup_{0\leq t\leq T}\left| \overline{f%
}_{m,n}^{i}(t,0)\right| \leq c,
\]
and
\[
\overline{\xi }^{m,n,1}\leq \overline{\xi }^{m,n,2},\mbox{ and }\overline{f}%
_{m,n}^{1}(t,y)\leq \overline{f}_{m,n}^{2}(t,y),\mbox{ for }(t,y)\in
[0,T]\times \mathbb{R,}
\]
in view of $\overline{f}^{1}(t,0)=f^{1}(t,0)=f^{2}(t,0)=\overline{f}%
^{2}(t,0) $. Using the comparison theorem \ref{compr4}-(i), we have for $%
0\leq s\leq t\leq T$%
\[
\overline{Y}_{t}^{m,n,1}\leq \overline{Y}_{t}^{m,n,2},\overline{K}%
_{t}^{m,n,1,+}-\overline{K}_{s}^{m,n,1,+}\geq \overline{K}_{t}^{m,n,2,+}-%
\overline{K}_{s}^{m,n,2,+}.
\]

By the convergence results in the step 4 and step 5 of the proof of theorem
\ref{exist22}, let $m\rightarrow \infty $, then $n\rightarrow \infty $, we
get for $0\leq s\leq t\leq T$
\[
\overline{Y}_{t}^{1}\leq \overline{Y}_{t}^{2},\overline{K}_{t}^{1,+}-%
\overline{K}_{s}^{1,+}\geq \overline{K}_{t}^{2,+}-\overline{K}_{s}^{2,+}.
\]
Especially, with $s=0$, it follows $\overline{K}_{t}^{1,+}\geq \overline{K}%
_{t}^{2,+}$. Similarly $\overline{K}_{t}^{1,-}\leq \overline{K}_{t}^{2,-}$.

(ii) and (iii) are from comparison theorem \ref{compr4} -(ii) and (iii),
with approximation as in (i), so we omit it. $\square $

\textbf{Acknowledgement } Author thanks professor Jean-Pierre
Lepeltier for his helps when author worked on this paper.

\end{document}